\pdfoutput=1
\documentclass[toc]{mathprint}
\usepackage{rutarmacros}
\usepackage{macros}
\title[Assouad Interpolation]{Interpolating with generalized Assouad dimensions}

\author[Banaji]
  {Amlan Banaji}
  {Department of Mathematics and Statistics, University of Jyväskylä, P.O.\ Box 35 (MaD), FI-40014 University of Jyväskylä, Finland}
  {banajimath@gmail.com}

\author[Rutar]
  {Alex Rutar}
  {Department of Mathematics and Statistics, University of Jyväskylä, P.O.\ Box 35 (MaD), FI-40014 University of Jyväskylä, Finland}
  {alex@rutar.org}

\author[Troscheit]
  {Sascha Troscheit}
{Department of Mathematics, Uppsala University, Box 480, 751 06 Uppsala, Sweden}
  {sascha.troscheit@math.uu.se}

\begin{document}
\begin{abstract}
    The $\phi$-Assouad dimensions are a family of dimensions which interpolate between the upper box and Assouad dimensions.
    They are a generalization of the well-studied Assouad spectrum with a more general form of scale sensitivity that is often closely related to ``phase-transition'' phenomena in sets.

    In this article we establish a number of key properties of the $\phi$-Assouad dimensions which help to clarify their behaviour.
    We prove for any bounded doubling metric space $F$ and $\alpha\in\R$ satisfying $\dimuB F<\alpha\leq\dimA F$ that there is a function $\phi$ so that the $\phi$-Assouad dimension of $F$ is equal to $\alpha$.
    We further show that the ``upper'' variant of the dimension is fully determined by the $\phi$-Assouad dimension, and that homogeneous Moran sets are in a certain sense generic for these dimensions.

    Further, we study explicit examples of sets where the Assouad spectrum does not reach the Assouad dimension.
    We prove a precise formula for the $\phi$-Assouad dimensions for the boundary of Galton--Watson trees that correspond to a general class of stochastically self-similar sets, including Mandelbrot percolation.
    The proof of this result combines a sharp large deviations theorem for Galton--Watson processes with bounded offspring distribution and a general Borel--Cantelli-type lemma for infinite structures in random trees.
    Finally, we obtain results on the $\phi$-Assouad dimensions of overlapping self-similar sets and decreasing sequences with decreasing gaps.
\end{abstract}

\section{Introduction}
A common theme in geometric measure theory and fractal geometry is to understand the fine scaling properties of sets in the absence of a smooth or differentiable structure.
An important quantity in this context is the notion of Assouad dimension.
This definition of dimension was first explicitly introduced by Assouad \cite{zbl:0396.46035} in his study of bi-Lipschitz embeddings of general metric spaces into Euclidean space.
The Assouad dimension also appears naturally as the maximal Hausdorff dimension of limits given by ``zooming in'' on the set; this work goes back to the pioneering work of Furstenberg and his notion of star dimension.
Especially in the past few decades, the Assouad dimension has received widespread attention from various perspectives: we refer the reader to the books on fractal geometry \cite{Fraser2021InterpolationSurvey}, embedding theory \cite{zbl:1222.37004}, and quasiconformal geometry \cite{zbl:1201.30002} for more details and background on this subject.

To make our context precise, we work with a non-empty, totally bounded metric space $(F,d)$.
Given $E\subseteq F$, we denote by $N_r(E)$ the smallest number of open balls of radius $r$ required to cover $E$.
Now, the \emph{Assouad dimension} of the set $F$ is given by
\begin{equation}
    \begin{aligned}\label{e:Assouad}
        \dimA F=\inf\Bigl\{s:(\exists C>0)&(\forall 0<r\leq R<1)\\*
                                                  &\sup_{x\in F}N_r(F\cap B(x,R))\leq C \Bigl(\frac{R}{r}\Bigr)^s\Bigr\}.
    \end{aligned}
\end{equation}
We always assume that $(F,d)$ is doubling, or equivalently that $\dimA F < \infty$.
If $F$ is Ahlfors--David $s$-regular, then $\dimuB F=\dimA F=s$.
However, for many important classes of sets (for example, self-similar sets with overlaps, self-affine sets, and random sets), it can happen that $\dimuB F<\dimA F$.
In this situation, we know that at some resolutions and locations, the set $F$ will have ``larger than average'' scaling.
A natural question, and one which is often important in applications, is to understand at which resolutions this scaling occurs.
Answering this question precisely has played a key role in answering questions as disparate as Hölder distortion estimates \cite{zbl:1390.28019}, conformal dimension \cite{zbmath:07731255}, and $L^p$-improving properties of maximal operators and fractal local smoothing estimates \cite{zbmath:07732556,zbl:1526.42033,arxiv:2501.12805}.

\subsection{Generalized Assouad dimensions}\label{ss:gen-Assouad}
In this article, we study the question of the fine scaling properties of Assouad dimension in a general setting and for some important explicit families of sets.
Perhaps the first systematic approach to this problem was the introduction of the \defn{Assouad spectrum} by Fraser \& Yu \cite{zbl:1390.28019}.
This modification of the Assouad dimension imposes the relationship $r=R^{1/\theta}$ for some fixed $\theta\in(0,1)$ in the definition of the Assouad dimension and is part of a more general scheme of ``dimension interpolation'' \cite{Fraser2021InterpolationSurvey}\footnote{Other notable examples include the \emph{(generalized) intermediate dimensions} \cite{zbl:1448.28009,zbmath:07750829} and the \emph{Fourier dimension spectrum} \cite{zbmath:07932394}.}.
The Assouad spectrum, which we denote by $\dimAs\theta F$, is a continuously parameterized family of
dimensions with $\dimuB F\leq\dimAs\theta F\leq\dimA F$ and satisfies $\lim_{\theta\to 0}\dimAs\theta F=\dimuB F$.

One might hope that the Assouad spectrum provides a robust theory of ``interpolation''; however, it need not hold that $\lim_{\theta\to 1}\dimAs\theta F=\dimA F$.
In fact, the quantity $\lim_{\theta\to 1}\dimAs\theta F$ coincides with the \emph{quasi-Assouad dimension}, denoted $\dimqA F$ and introduced in \cite{zbl:1345.28019}.
It is possible for the quasi-Assouad dimension to be smaller than the Assouad dimension for two well-studied classes of sets:
\begin{enumerate}[itemsep=1mm]
    \item Random sets, which occur naturally as limiting objects resulting from branching processes (see \cref{s:stoch}).
    \item Dynamically invariant sets exhibiting some form of overlapping behaviour, such as the invariant sets of overlapping self-similar iterated function systems in $\R$ (see \cref{s:other-ex}).
\end{enumerate}
Within these families, it appears that the ``expected'' behaviour is that $\dimAs\theta F=\dimuB F$ for all $\theta\in(0,1)$, but $\dimA F$ is as large as possible.
(It is an interesting, and probably challenging, open question to verify if $\dimAs\theta F=\dimuB F$ for all $\theta\in(0,1)$ for \emph{all} self-similar sets $F$.)
Moreover, for self-affine sets $F \subset \R^2$, it can happen that $\dimqA F<\dimA F$ even in the strongly separated case, see~\cite{zbmath:07808129}.

As a way to remedy this situation, Fraser \& Yu suggested in \cite{zbl:1390.28019} that one might generalize the Assouad spectrum by instead allowing the smaller scale $r$ to be prescribed as a function of the larger scale $R$, for some sensibly-chosen but otherwise arbitrary function $\phi$.
This program was taken up by Garc\'ia, Hare, \& Mendivil in \cite{zbl:1485.28006}, who established various fundamental properties of this construction.
We also refer the reader to \cite{zbl:1437.28015, zbl:1497.28005} for other articles studying the generalized Assouad dimensions of some specific families of sets and to \cite[§3.3.3]{Fraser2021InterpolationSurvey} for more introduction to this program\footnote{A version of the generalized Assouad dimensions for measures has been studied in \cite{zbl:1506.28003,zbl:1489.28001,zbl:1543.28003}.}.

We define a particular variant here, which is trivially more restrictive than the original definition but for all purposes functions in the same way.
\begin{definition}\label{d:dimf}
    We say that a function $\phi\colon (0,1)\to\R^+$ is a \defn{dimension function} if the following two conditions hold:
    \begin{enumerate}[nl,r]
        \item\label{im:phi-gap} $R\mapsto\phi(R)\log(1/R)$ increases to infinity as $R$ decreases to zero, and
        \item\label{im:phi-decr} $\phi(R)$ decreases as $R$ decreases to zero.
    \end{enumerate}
\end{definition}
Note that these conditions necessarily imply continuity of $\phi$.
We denote the set of all dimension functions by $\mathcal{W}$.
For $\phi\in\mathcal{W}$, we then define the \defn{$\phi$-Assouad dimension}\footnote{It seems that every sensible choice of specification (e.g.~$r=\phi(R)$, $r=R^{\phi(R)}$, or $r=R^{1+\phi(R)}$) has its own share of benefits and drawbacks.
One must be careful when reading the literature to verify which notation is used.}
of $F$ by
\begin{align*}
    \dimAs\phi F=\inf\Bigl\{s:(\exists C>0)&(\forall 0<r=R^{1+\phi(R)}\leq R<1)\\*
                                           &\sup_{x\in F}N_r(F\cap B(x,R))\leq C R^{-\phi(R) s}\Bigr\}.
\end{align*}
We similarly define the \defn{upper $\phi$-Assouad dimension}\footnote{
What we call the upper $\phi$-Assouad dimension is referred to as the $\phi$-dimension in \cite{zbl:1485.28006}. Our terminology is chosen for consistency with the (upper) Assouad spectrum; see \cite{zbl:1410.28008}.}
of $F$ by
\begin{align*}
    \dimuAs\phi F=\inf\Bigl\{s:(\exists C>0)&(\forall 0<r\leq R^{1+\phi(R)}\leq R<1)\\*
                                           &\sup_{x\in F}N_r(F\cap B(x,R))\leq C \Bigl(\frac{R}{r}\Bigr)^s\Bigr\}.
\end{align*}
Occasionally, we will write $\dimAs\phi F$ to denote the identical formula whenever $\phi\colon(0,1)\to\R^+$ is any function.
One convenient feature of the $\phi$-Assouad dimension is that it can be expressed as a limit: it follows directly from the definition that
\begin{equation}\label{e:phi-dimA-lim}
    \dimAs\phi F=\limsup_{R\to 0}\frac{\log\sup_{x\in F}N_{R^{1+\phi(R)}}\bigl(F\cap B(x,R)\bigr)}{-\phi(R)\log R}.
\end{equation}

Intuitively, condition \cref{im:phi-gap} of $\phi$ being a dimension function means that the gap between the upper and lower scales grows monotonically and arbitrarily large as the scale goes to zero, and condition \cref{im:phi-decr} means that the $\phi$-Assouad dimension at scale $R$ becomes ``more like the Assouad dimension'' as $R$ tends to $0$.
This analogy is made more precise in \cref{ss:interpolation-recover}.
If $\phi$ is a function such that $\phi(R)\log(1/R)$ is increasing but does not diverge to infinity, then we would have $\dimAs\phi F=0$ for all bounded sets $F$, whereas we want to insist that $\dimAs\phi F\geq\dimuB F$.
Dimension functions are abundant: it is proven in \cref{p:dim-maximal} that for any function $\phi$ satisfying \cref{im:phi-gap}, there is a unique maximal dimension function $\psi\leq\phi$ (with the partial order of pointwise comparison).

Our main contributions in this article are four-fold\footnote{For clarity of exposition, we will not consider the analogous results for the dual notion of $\phi$-lower dimensions, though we expect many of our proofs to work in a similar way in that situation.
More detail on the $\phi$-lower dimensions can be found in \cite{zbl:1485.28006}.}:
\begin{enumerate}[itemsep=1mm]
    \item We establish and clarify general properties of the $\phi$-Assouad dimensions.
    \item We show that the $\phi$-Assouad dimensions recover the interpolation: for any $\alpha\in\R$
      satisfying $\dimuB F<\alpha\leq\dimA F$, there is a dimension function $\phi$ so that $\dimAs\phi F=\alpha$.
    \item We establish precise formulas for the $\phi$-Assouad dimensions of stochastically self-similar sets (in particular, the Gromov boundary of Galton--Watson processes with finite support).
        These results are consequences of more general sharp results on large deviations of Galton--Watson processes, which may be of independent interest.
    \item We investigate general properties of overlapping self-similar sets, and prove initial quantitative bounds on the $\phi$-Assouad dimensions for some specific examples.
\end{enumerate}
We will discuss our main results along these themes more precisely for the remainder of this introduction.

\subsection{Rate windows and regularity of \texorpdfstring{$\phi$}{𝜙}-Assouad dimension}\label{ss:windows}
A particularly important concept in the notion of the $\phi$-Assouad dimension is the definition of a \defn{dimension rate window}.
As we will see in \cref{ex:spectrum-equiv} below, this notion generalizes the definition of the Assouad spectrum.
Within the rate window corresponding to a dimension function, the corresponding set of dimensions is relatively well-behaved.
Moreover, one might hope that for sufficiently nice sets, there is exactly one dimension rate window on which the $\phi$-Assouad dimensions exhibit non-endpoint behaviour.
As we will see, this is in fact the case for the $\phi$-Assouad dimensions of stochastically self-similar sets.

In this section, we discuss some general ideas behind the notion of the rate window; we hope that this will help to clarify the statement of the results in the following section.
Fix a dimension function $\psi\in\mathcal{W}$.
For $\alpha\in(0,\infty)$, we denote by $\psi_\alpha$ the function $x\mapsto\psi(x)/\alpha$.
Observe that $\psi_\alpha\in\mathcal{W}$, and $R^{1+\psi_\alpha(R)}$ increases as $\alpha$ increases.
We then define the \defn{dimension rate window} of $\psi$ by
\begin{equation*}
    \mathcal{W}_\psi=\bigl\{\psi_\alpha:\alpha\in(0,\infty)\bigr\}\subset\mathcal{W}
\end{equation*}
To recall a familiar example, consider the case when $\psi$ is a constant function.
This corresponds to the usual Assouad spectrum.
\begin{example}\label{ex:spectrum-equiv}
    Suppose $\psi(x)=\frac{1}{\theta}-1$, so that $\dimAs\psi$ is precisely the usual Assouad spectrum at $\theta$.
    Then
    \begin{equation*}
        W_\psi=\bigl\{x\mapsto\frac{1}{\alpha}\Bigl(\frac{1}{\theta}-1\Bigr):\alpha\in(0,\infty)\bigr\}.
    \end{equation*}
    In particular, suppose $\theta_\alpha$ is chosen so that
    \begin{equation}\label{e:Assouad-rel}
        \frac{1}{\theta_\alpha}-1=\frac{1}{\alpha}\Bigl(\frac{1}{\theta}-1\Bigr).
    \end{equation}
    A direct computation shows that $\theta_\alpha$ is an increasing function of $\alpha$ with $\lim_{\alpha\to 0}\theta_\alpha=0$, $\theta_1=\theta$, and $\lim_{\alpha\to\infty}\theta_\alpha=1$.
    In other words, $\mathcal{W}_\psi$ is simply the set of dimension functions corresponding to the usual Assouad spectrum, for all $\theta\in(0,1)$.
\end{example}
As detailed in the fundamental result \cref{it:fund-dim}, many of the properties that one expects for the Assouad spectrum generalize to the setting of arbitrary dimension rate windows.

Recall that a metric space $F$ is \emph{doubling} if the number of $r/2$-balls needed to cover any $r$-ball centred in $F$ is bounded above by a uniform integer $M_F$; the smallest such $M_F$ is called the \emph{doubling constant}.
We work with a non-empty, bounded, doubling (abbreviated n.b.d.) metric space $F$.
These conditions imply that $F$ is totally bounded (so $N_r(B(x,R))$ is always finite) and has finite Assouad dimension.
For some results we will specialize to the case when $F$ is a non-empty, bounded subset of $\R^d$.
\begin{itheorem}\label{it:fund-dim}
    Let $\phi$ and $\psi$ be dimension functions and let $F$ be any n.b.d.\ space.
    \begin{enumerate}[nl,r]
        \item\label{im:limit-eq} If $\displaystyle\lim_{R\to 0}\frac{\phi(R)}{\psi(R)}=\alpha\in(0,\infty)$, then $\dimAs\psi F=\dimAs{\phi_\alpha} F$.
        \item \label{im:limit-diff} If $\dimAs\phi M = \dimAs\psi M$ for all non-empty, compact, perfect sets $M\subset\R$ then $\lim_{R\to 0}\frac{\phi(R)}{\psi(R)}=1$.
        \item\label{im:window-compare} If $\displaystyle\lim_{R\to 0}\frac{\phi(R)}{\psi(R)}=0$, then $\dimuAs\psi F\leq\dimAs\phi F$.
        \item\label{im:cont} $\alpha\mapsto\dimAs{\phi_\alpha}F$ is a continuous function of $\alpha$.
    \end{enumerate}
\end{itheorem}
The proof of this result is given in \cref{ss:window-proofs}.
Note that \cref{im:limit-eq} essentially says that we may have instead defined
\begin{equation*}
    \mathcal{W}_\psi=\bigl\{\phi\in\mathcal{W}:\lim_{R\to 0}\frac{\phi(R)}{\psi(R)}\in(0,\infty)\bigr\}.
\end{equation*}
Facts \cref{im:limit-eq} and \cref{im:cont} were originally proven in \cite[Section~1.3.1,~(3)]{zbl:1485.28006} and \cite[Section~1.3.1,~(7)]{zbl:1485.28006} respectively for the \emph{upper} $\phi$-Assouad dimensions.
Our proofs of \cref{im:limit-eq} (given in \cref{c:phi-equivs}) and \cref{im:cont} (given in \cref{p:window-bounds}) follow similarly, with the only additional complication being non-monotonicity of $\alpha\mapsto\dimAs{\phi_\alpha}F$.

Statement \cref{im:limit-diff} is similar to \cite[Theorem~3.8]{zbl:1485.28006}, though since we use the $\phi$-Assouad dimensions which precisely specify the relationship between the scales, we obtain a converse for \cref{im:limit-eq}.
The proof of \cref{im:limit-diff} is given in \cref{p:continuity-converse} and \cref{c:phi-equivs} by a direct argument using a Moran construction.

Finally, \cref{im:window-compare} is proved in \cref{p:inter-window-compare}, and does not seem to have been observed before.
Heuristically, \cref{im:window-compare} states that dimension functions in distinct windows yield notions of dimension which satisfy a strong ordering property.
Our \cref{it:fund-dim} also gives new intuition for the observation that the upper box dimension is a lower bound for the $\phi$-Assouad dimension (and in particular, the Assouad spectrum), and why allowing scales $r$ and $R$ close together in \cref{e:Assouad} increases the corresponding dimensional constant: repeating good bounds for close scales yields good bounds for well-separated scales.
The details are given in \cref{t:Assouad-recover}.

Another application of \cref{it:fund-dim} is that it gives a certain way to understand the space of $\psi$-Assouad dimensions.
Define an equivalence relation $\sim$ on the space of dimension functions by $\phi \sim \psi$ if $\lim_{R\to 0}\phi(R)/\psi(R) =1$.
By \cref{c:phi-equivs}, the set $\mathcal{D}$ of equivalence classes precisely corresponds to the set of different notions of $\psi$-Assouad dimension.
We can define a natural non-strict partial order $\preceq$ on $\mathcal{D}$ by $[\phi] \preceq [\psi]$ if $\dimAs\phi F \leq \dimAs\psi F$ for all bounded $F \subset \R^d$.
By \cref{it:fund-dim} \cref{im:limit-diff} (or more precisely \cref{p:continuity-converse}) and \cref{im:cont}, $[\phi] \preceq [\psi]$ if and only if $\limsup_{R \to 0} \psi(R) / \phi(R) \leq~1$.
The natural topology $\mathcal{T}$ on $\mathcal{D}$ is the initial topology of the set of functions
\begin{equation*}
    \bigl\{  f_F : d \in \N, F \subset \R^d \mbox{ bounded} \bigr\},
\end{equation*}
where $f_F \colon \mathcal{D} \to \R$, $f_F([\psi]) = \dimAs\psi F$, and where we endow $\R$ with its usual topology.
We have the following explicit description of $\mathcal{T}$, which is proven at the end of \cref{ss:window-proofs}.
\begin{icorollary}\label{ic:topology}
    A basis of open sets for the topology $\mathcal{T}$ is
    \begin{equation*}
        \{  N_{\psi,\eps} : \psi\in\mathcal{W}, \eps \in (0,1)  \},
    \end{equation*}
    where
    \begin{equation*}
        N_{\psi,\eps} \coloneqq \left\{  [\phi] \in \mathcal{D} : 1-\eps < \liminf_{R \to 0} \frac{\phi(R)}{\psi(R)} \leq \limsup_{R \to 0} \frac{\phi(R)}{\psi(R)} < 1+\eps  \right\}.
    \end{equation*}
\end{icorollary}
One can define an equivalence relation $\sim_w$ on $\mathcal{D}$ by $[\phi] \sim_w [\psi]$ if
\begin{equation*}
    0 < \liminf_{R \to 0} \frac{\phi(R)}{\psi(R)} \leq  \limsup_{R \to 0} \frac{\phi(R)}{\psi(R)} < \infty.
\end{equation*}
The equivalence classes of $\sim_w$ are clearly elements of $\mathcal{T}$, and they partition $\mathcal{D}$.
Moreover, there are uncountably many such equivalence classes, for example the equivalence classes of each $R\mapsto(-\log R)^{-t}$ for $0 < t < 1$ (this example is discussed in \cite[Theorem~3.9]{zbl:1485.28006}).
In particular, $(\mathcal{D},\mathcal{T})$ is neither connected nor separable.

\subsection{Main results for general sets}\label{ss:intro-general}
Now that we have introduced our setting and established some general notation, we begin by stating our main results.
Our first two results are general facts about the $\phi$-Assouad dimensions.

First, the upper $\phi$-Assouad dimensions can be obtained from the $\phi$-Assouad dimensions, which is a generalization of the corresponding result for the Assouad spectrum, \cite[Theorem~2.1]{zbl:1410.28008}.
We recall in general that the $\phi$-Assouad dimension and the upper $\phi$-Assouad dimension need not be equal; see, for instance, \cite[§8]{zbl:1390.28019}.
\begin{itheorem}\label{it:upper-recover}
    Let $F$ be an n.b.d.\ space and let $\phi$ be a dimension function.
    Then
    \begin{equation*}
        \dimuAs\phi F=\sup_{\alpha\in(0,1)}\dimAs{\phi_\alpha} F.
    \end{equation*}
\end{itheorem}
The proof of this result is given in \cref{ss:upper-recover} and uses \cref{it:fund-dim}~\cref{im:window-compare} in a critical way.

Next, we show that the $\phi$-Assouad dimensions \emph{recover the interpolation}.
\begin{itheorem}\label{it:interpolation-recover}
    For any n.b.d.\ space $F$ and $\dimuB F < \alpha \leq \dimA F$, there exists a dimension function $\phi$ such that
    \begin{equation*}
        \dimuAs\phi F=\dimAs\phi F=\alpha.
    \end{equation*}
\end{itheorem}
The proof of this result is given in \cref{ss:interpolation-recover} and consists of the main technical work of our results on general sets.
In fact, as detailed in \cref{t:Assouad-recover}, for the specific value $\alpha=\dimA F$ we can choose the dimension function $\phi$ so that it is essentially as small as reasonably possible.
This implies that one may take the pairs of scales to be arbitrarily close together in the definition of the Assouad dimension, as long as the ratio tends to zero.
Note that the $\phi$-Assouad dimensions were one of the motivations for the introduction of the generalized intermediate dimensions, and some of our results have parallels in that setting.
In particular, \cite[Theorem~6.1]{zbmath:07750829} says that the generalized intermediate dimensions can be used to recover the interpolation between Hausdorff and box dimension for any compact set, which is analogous to \cref{it:interpolation-recover}.

In \cref{ss:moran} we provide further results concerning the general construction of sets which realize the $\phi$-Assouad dimensions.
We defer the precise definition of a homogeneous Moran set to \cref{ss:window-proofs}: heuristically, these sets have maximal homogeneity in space (in other words, the set looks the same everywhere) but need not have any homogeneity between scales.
A canonical example is a Cantor set in $\R$ formed by contracting with ratio $\alpha\in(0,1/2]$ at some steps, and with ratio $\beta\in(0,1/2]$ at the other steps.
Homogeneous Moran sets were used in \cite[Theorem~3.9]{zbl:1485.28006} to give the first example of a set $M$ with $\dimqA M < \dimA M$ for which for all $s \in [\dimqA M, \dimA M]$ there exists $\phi_s$ with $\dimAs{\phi_s} M = s$.
In \cref{p:moran-dim-formula} we give a short proof of a formula for the $\phi$-Assouad dimensions of homogeneous Moran sets, which may be of interest in its own right.

In \cite{zbmath:07937992}, homogeneous Moran sets are used to characterize the attainable forms of Assouad spectra.\footnote{A possible direction for future research would be to try to characterize the attainable forms of $\phi$-Assouad dimensions within a given window, but we will not pursue this.}
The techniques are based on constructions developed for characterizing the intermediate dimensions, which appeared in \cite{zbl:1509.28005}.
A corollary of the results therein is that for any bounded set $F$, there is a homogeneous Moran set $M$ such that for all $\theta\in(0,1)$, $\dimAs\theta F=\dimAs\theta M$.
It is natural to ask to what extent the behaviour of the $\phi$-Assouad dimensions for different functions $\phi$ can be witnessed by homogeneous Moran sets.
In particular, we ask the following question which we are unable to answer in full generality.
\begin{question}
    Given $d \in \N$ and $F \subset \R^d$, does there necessarily exist a homogeneous Moran set $M \subset \R^d$ such that $\dimAs{\phi} F = \dimAs{\phi} M$ for all dimension functions $\phi$?
\end{question}
We make good progress towards an affirmative answer by proving the following result in \cref{ss:moran}, which says that homogeneous Moran sets are indeed typical for families of dimension functions which satisfy mild conditions.
These conditions hold for many large families; for example, the ordering condition on each $\mathcal{A}_i$ holds for the family used in the proof of \cite[Theorem~3.9]{zbl:1485.28006}.
A particular example of a family for which it holds is $\{ \phi_t(R) = (\log (1/R))^t\}_{0 < t < 1}$.
\begin{itheorem}\label{it:morantypical}
    Fix $d \in \N$ and $F \subset \R^d$, and let $\mathcal{A}$ be a family of dimension functions.
    Suppose $\mathcal{A}=\bigcup_{i=1}^\infty\mathcal{A}_i$ where for each $i$ there exists $T_i \subset \R$ such that $\mathcal{A}_i = \{ \phi_{i,t} : t \in T_i\}$ and whenever $t,t' \in T_i$ satisfy $t \geq t'$ the following limit exists and lies in $[0,1]$:
    \begin{equation}
        \lim_{R\to 0}\frac{\phi_{i,t}(R)}{\phi_{i,t'}(R)}\in[0,1].
    \end{equation}
    Then there exists a homogeneous Moran set $M \subset \R^d$ such that $\dimAs\psi M \leq\dimAs\psi F$ for all dimension functions $\psi$, and moreover
    \begin{equation*}
        \dimAs{\psi} F = \dimAs{\psi} M \qquad \text{and} \qquad \dimuAs{\psi} F = \dimuAs{\psi} M
    \end{equation*}
    for all $\phi \in \mathcal{A}$ and $\psi\in\mathcal{W}_{\phi}$.
\end{itheorem}
\begin{remark}
    Taking some $\phi\in \mathcal{A}$ to be constant, we can guarantee that $\dimuB F = \dimuB M$, $\dimqA F = \dimqA M$, and $\dimAs{\theta} F = \dimAs{\theta} M$ for all $\theta \in (0,1)$.
    In particular, \cref{it:morantypical} gives a direct proof of the fact from \cite{zbmath:07937992} that all possible behaviours of Assouad spectra can be realized by homogeneous Moran sets.
\end{remark}

\subsection{Main results for specific sets}\label{ss:intro-specific}
We now turn our attention to some specific families of sets.
As discussed earlier, our primary motivation for studying the $\phi$-Assouad dimensions is that many
natural families of sets exhibit the dichotomy that the Assouad dimension is as large as possible, whereas the Assouad spectrum is constantly equal to the box dimension.

\subsubsection{Branching processes}
Our first consideration, and the situation in which we have the most precise and general results, is on branching processes and associated random sets.
We refer the reader to \cref{ss:gw-proc} for more precise definitions and only give an overview in this section.
Let $X$ be a non-negative integer valued random variable; assume that $X$ has finite support.
Then the Galton--Watson process $Z_k$ with offspring variable $X$ is defined by the recursion
\begin{equation*}
    Z_0=1\qquad\text{and}\qquad Z_{k+1} = \sum_{i=1}^{Z_k} X_{k,i}
\end{equation*}
where $X_{k,i}$ are independent random variables with the same distribution as $X$.
In other words, at step $k+1$, every child which appeared at step $k$ yields a random number of offspring controlled by the random variable $X$.
Assuming that $\E(X)>1$, iterating this process yields an infinite tree associated with the Galton--Watson process, and this tree has a \defn{Gromov boundary} $\partial\cT$ that is non-empty with positive probability.
Fixing the natural metric from the longest common substring makes the Gromov boundary a doubling metric space.

Note that, almost surely conditioned on non-extinction, the tree $\partial\cT$ will have a subtree with full branching over $k$ steps for any $k\in\N$ with a surviving branch in every such leaf.
In particular, this implies that the Assouad dimension will be as large as possible.
Basic scaling arguments further show that the box-counting and Hausdorff dimension coincide and equal the ``average'' exponential growth rate, see e.g.\ \cite{zbl:1437.28015}.

The question of determining the $\phi$-Assouad dimension then (morally) reduces to the following question: at which length scales do we begin to see larger than expected subtrees?
Key to this problem is the large deviations of the underlying branching process.
In particular, we use the following large deviations theorem, which is more refined than the Chernoff type bound in \cite[Lemma~3.1]{zbl:1437.28015}.
Here and elsewhere, we write $A \lesssim B$ to mean $A \leq CB$ for some implicit constant $C$, and if $A \lesssim B$ and $B \lesssim A$ we write $A \approx B$.
\begin{proposition}\label{p:GW-large-deviations}
    Let $Z_k$ be a Galton--Watson process with offspring random variable $X$ which is not almost surely constant.
    Assume that its probability generating function $f$ is a polynomial of degree $2\leq N <\infty$ and $m\coloneqq \E(X)>1$.
    Define $\gamma$ such that $m^\gamma = N$.
    Then for all $1<t<\gamma$, all $\eps>0$ sufficiently small, and all $k \in \N$,
    \begin{equation*}
        \exp\left( -m^{(t-1+\eps)\frac{\gamma}{\gamma-1}k} \right)\lesssim
        \Prob\left( Z_k  \geq m^{tk} \right)
        \lesssim \exp\left( -m^{(t-1-\eps)\frac{\gamma}{\gamma-1}k} \right),
    \end{equation*}
    with the implicit constants depending only on $t$ and $\eps$.
\end{proposition}
Large deviation results of this type have been known since the seminal work of Biggins \& Bingham \cite{zbl:0796.60090} and we discuss its history and variants in \cref{s:large-deviations}.
However, in order to obtain our desired results on the $\phi$-Assouad dimension of random trees, we need something somewhat different than a large deviations estimate: we must guarantee almost sure occurrence of infinitely many subtrees which are smaller (or larger) than expected.
The key ingredient here is the following Borel--Cantelli lemma for trees.
\begin{lemma}\label{l:tree-bc}
    Let $E_k$ be any measurable event for a Galton--Watson tree and write $P_k = \Prob(E_k)$.
    Let $\widetilde{E}$ be the event that there are infinitely many $k\in\N$ such that a Galton--Watson
    tree contains a subtree $\cT(v)\in E_k$ at level $k$.
    \begin{enumerate}[nl,r]
        \item $\Prob(\widetilde{E}) =0$ \;if\; $\sum_{n\in\N}P_n m^n <\infty$,
        \item $\Prob(\widetilde{E}) =1$, conditioned on non-extinction, if there exists a summable sequence $K_n$ of non-negative numbers such that $\sum_{n\in\N}K_n P_n m^n=\infty$.
    \end{enumerate}
\end{lemma}
By combining \cref{p:GW-large-deviations} and \cref{l:tree-bc}, we obtain the following sharp result for the $\phi$-Assouad dimensions.
\begin{itheorem}\label{it:tree-dims}
    Let $Z_k$ be a Galton--Watson process with finitely supported offspring distribution with mean $m$ and maximal offspring number $N$.
    Let $\partial\cT$ denote the Gromov boundary of the associated Galton--Watson tree.
    Write
    \begin{equation}\label{e:psi-def}
        \psi(R)=\frac{\log\log(1/R)}{\log(1/R)}.
    \end{equation}
    The following results hold almost surely conditioned on non-extinction.

    For any dimension function $\phi$, if $\lim_{R\to 0}\frac{\psi(R)}{\phi(R)}=\alpha\in[0,\log N]$, then
    \begin{equation*}
        \dimAs\phi \partial\cT = \alpha\left(1-\frac{\log m}{\log N}\right) +\log m.
    \end{equation*}
    Otherwise, if $\lim_{R\to 0}\frac{\psi(R)}{\phi(R)}\geq\log N$, then
    \begin{equation*}
        \dimAs\phi\partial\cT =\log N.
    \end{equation*}
\end{itheorem}
We recall that the almost sure box and Assouad dimensions are given by $\log m$ and $\log N$, respectively.
In the terminology introduced in \cref{ss:windows}, the dimension rate window defined by the function $\psi$ in \cref{e:psi-def} fully recovers the interpolation between the box and Assouad dimensions.

\begin{figure}[t]
    \centering
    \begin{subcaptionblock}{.47\textwidth}
        \centering
        \includegraphics[width = 0.9\textwidth]{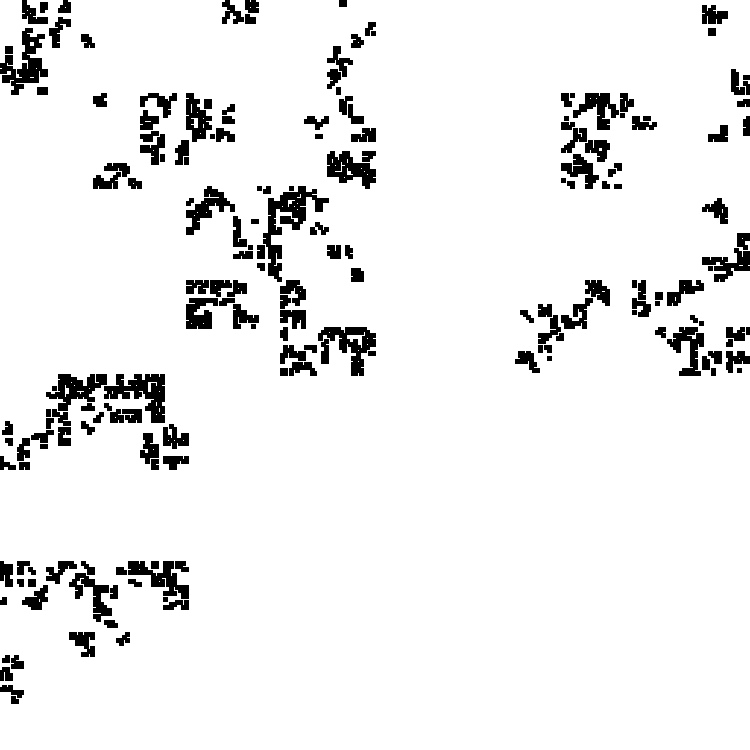}
        \caption{Unconditioned.}
    \end{subcaptionblock}
    \begin{subcaptionblock}{.47\textwidth}
        \centering
        \includegraphics[width = 0.9\textwidth]{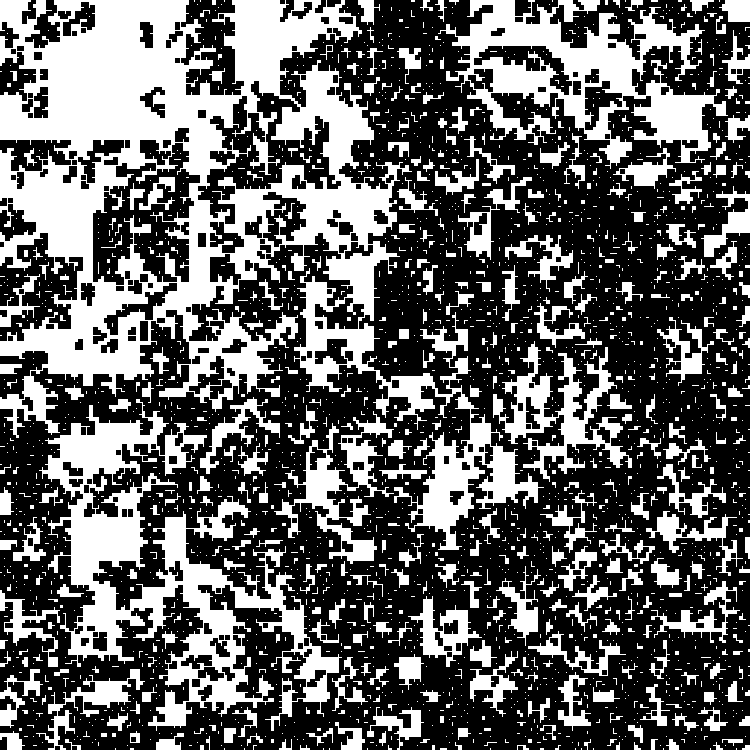}
        \caption{Conditioned.}
        \label{f:conditioned}
    \end{subcaptionblock}
    \caption{Typical Mandelbrot percolation of the unit square with parameters $n=2$, $d=2$, $p=0.65$.
    (B) is conditioned on every subcube surviving the first $6$ iteration steps.}
    \label{fig:Mandelbrot}
\end{figure}

In particular, we can apply \cref{it:tree-dims} to Mandelbrot percolation of the $d$-dimensional unit cube.
This is the limit set $M$ obtained by subdividing the unit cube in $\R^d$ into $n^d$ subcubes of side-length $n^{-1}$ and retaining each subcube independently with probability $p$, and then continuing the subdivision and retention ad infinitum.
Some instances of Mandelbrot percolation are depicted in \cref{fig:Mandelbrot}.
This is a well studied class of random fractals with a long history, see \cite{Rams2014mandelbrot,zbl:0289.76031} and \cite[Section~15.2]{zbl:1285.28011}.
More recently, the Assouad dimension has been studied in \cite{zbl:1392.37024} building on \cite{zbl:1411.28006}, and the Assouad spectrum has been studied in \cite{zbl:1390.28019,zbl:1437.28015}.
Almost surely, $\dimH M = \dimB M = \dimAs\theta M = d + \log p / \log n$ while $\dimA M = d$.

Fix the function $\psi$ defined in \cref{e:psi-def}.
\begin{icorollary}\label{ic:percolation}
    Let $M$ be a Mandelbrot percolation set with retention probability $p\in (0,1]$, division parameter $n\geq 2$ and ambient space dimension $d$ such that $pn^d>1$.
    The following results hold almost surely conditioned on non-extinction.
    For any dimension function $\phi$, if $\lim_{R\to 0}\frac{\psi(R)}{\phi(R)}=\alpha\in[0,\log n^d]$,
    \begin{equation*}
        \dimAs\phi M = \alpha\frac{\log (1/p)}{d\log^2 n} +\frac{\log pn^d}{\log n}.
    \end{equation*}
    Otherwise, if $\liminf_{R\to 0}\frac{\psi(R)}{\phi(R)}\geq\log n^d$, then
    \begin{equation*}
        \dimAs\phi M = d.
    \end{equation*}
\end{icorollary}
This gives a complete answer to \cite[Questions~17.6.1 \& 17.6.2]{Fraser2021InterpolationSurvey} and closes the gap in \cite[Theorem~2.3]{zbl:1437.28015}.

Interestingly, the large deviation formula in \cref{p:GW-large-deviations} and the dimension results in \cref{it:tree-dims} only depend on the maximal number of offspring and the mean offspring number, rather than the entire offspring distribution.
For example, fix parameters in Mandelbrot percolation so that $pn^d>1$.
If we now allow different hypercubes to have different retention probabilities, as long as we do this in a way that keeps all of the probabilities strictly positive and leaves their sum unchanged at $pn^d$, the dimension formulas will almost surely stay the same.
This is the case even if the product of the probabilities (i.e.\ the probability of retaining all $n^d$ hypercubes at a given level) changes.

\subsubsection{Self-similar sets}
Fix a finite index set $\mathcal{I}$ and for each $i\in\mathcal{I}$ a contracting similarity $S_i\colon\R\to\R$, i.e.\ $S_i(x)=r_i x+d_i$ where $0<|r_i|<1$.
Then the corresponding \defn{self-similar set} is the unique non-empty compact set $K$ satisfying
\begin{equation*}
    K=\bigcup_{i\in\mathcal{I}}S_i(K).
\end{equation*}
In \cite[Theorem~1.3]{zbl:1317.28014}, it is proven that if $K$ does not satisfy the \defn{weak separation condition} of Lau \& Ngai \cite{zbl:0929.28007}, then $\dimA K=1$.
The precise definition of the weak separation condition can be found in \cref{ss:similar-upper}; we just note here that if $\dimH K<1$, then the weak separation condition is equivalent to Ahlfors--David regularity which implies that $\dimH K = \dimA K$ \cite{zbl:1334.28012}.
In particular we are interested in the case that the weak separation condition fails.
The dimension theory of deterministic self-similar sets beyond the weak separation condition has been a core area of study in fractal geometry; for a (highly incomplete) list of some notable results, see \cite{zbl:1230.37031,zbl:1337.28015,zbl:1426.11079,zbl:1426.28024,zbmath:07835917,zbl:1408.28017}.

For notational simplicity in the introduction, suppose the IFS is \emph{equicontractive}, so that $r_i=r\in(0,1)$ for all $i\in\mathcal{I}$ (the general case is treated in \cref{ss:similar-upper}).
Then set
\begin{equation*}
    \widetilde{M}_n=\sup_{x\in K}\#\bigl\{ g:(i_1,\ldots,i_n)\in\mathcal{I}^n,g=S_{i_1}\circ\cdots\circ S_{i_n}, g(K)\cap B(x,r^n)\neq\varnothing \bigr\}.
\end{equation*}
The weak separation condition is equivalent to requiring that $\sup_{n\in\N}\widetilde{M}_n<\infty$.
By using the cylinders giving $\widetilde{M}_n$, we can prove a natural upper bound for the $\phi$-Assouad dimensions general self-similar sets.
We also exhibit a family of sets for which we can demonstrate a reasonable lower bound.
This is described in the following result:
\begin{itheorem}\label{it:similar-bounds}
    Let $K$ be an equicontractive self-similar set and suppose $\phi$ is any dimension function such that
    \begin{equation}\label{e:Mn-window}
        \lim_{n\to \infty}\frac{\log \widetilde{M}_n}{n\cdot\phi(r^n)}=0.
    \end{equation}
    Then $\dimAs\phi K=\dimH K$.

    On the other hand, for all $m\in\N$ with $m\geq 3$, there is an explicit self-similar set $K\subset\R$ with three maps, each with contraction ratio $1/m$, and a constant $C_m>0$, such that if $\phi$ is any dimension function satisfying
    \begin{equation*}
        \liminf_{n\to\infty}\frac{\log n}{n\cdot \phi\bigl(m^{-n}\bigr)}\geq C_m
    \end{equation*}
    then $\dimAs\phi K=\dimA K=1$.
\end{itheorem}
The proof can be obtained by combining \cref{p:self-similar-upper} and \cref{t:ss-lower}.
The main difficulty with proving good lower bounds is that the presence of a large number of cylinders (as governed by $\widetilde{M}_n$) is insufficient to guarantee large covering number locally.
In some sense, we require the stronger property that the cylinders corresponding to a maximal ball $B(x,r)$ for $\widetilde{M}_n$ are ``uniformly distributed'' so that each cylinder individually makes a maximal contribution.

Unfortunately, we are unable to give an explicit example of a set for which the upper bound and lower bounds give the same dimension window.
We also wonder if the natural upper bound defined by the $\widetilde{M}_n$ is sharp.
\begin{question}
    Suppose $K$ is an equicontractive self-similar set in $\R$ with $\dimH K<1$ which does not satisfy the weak separation condition.
    Does there necessarily exist a single dimension function $\psi$ such that
    \begin{equation*}
    \lim_{\alpha\to 0}\dimAs{\psi_\alpha}K=\dimH K \qquad \mbox{and} \qquad \lim_{\alpha\to\infty}\dimAs{\psi_\alpha}K=\dimA K?
    \end{equation*}
    If so, is the asymptotic threshold given by $\log\widetilde{M}_n/n$ (as in \cref{e:Mn-window})?
\end{question}
Of course, this holds automatically if $K$ satisfies the weak separation condition (which implies that $\dimH K=\dimA K$) or if $\dimH K$ is equal to the ambient dimension.

The set for which we construct a non-trivial lower bound for \cref{it:similar-bounds} satisfies the exponential separation condition introduced in \cite{zbl:1337.28015}.
In particular, combined with work in \cite{zbl:1426.11079}\footnote{This argument is given explicitly in \cite[Theorem~1.3]{zbl:1469.28003} and is implicit in the proof of \cite[Theorem~4.1]{zbl:1407.28002}.}, $\widetilde{M}_n$ must grow subexponentially, i.e.
\begin{equation*}
    \lim_{n\to\infty}\frac{\log \widetilde{M}_n}{n}=0.
\end{equation*}
However, in order to prove that the lower bound in \cref{it:similar-bounds} is sharp, we must in fact prove the stronger upper bound (for that specific example) that $\widetilde{M}_n$ grows at most polynomially in $n$, i.e.\ $\widetilde{M}_n\leq C n^t$ for some $t>0$ and some $C>0$.

\subsubsection{Decreasing sequences with decreasing gaps}
Finally, in \cref{ss:decr-seq} we prove a sharp result describing the $\phi$-Assouad dimension of decreasing sequences with decreasing gaps.
We refer the reader to that section for the precise formulation and statement of the result.

\section{General properties of \texorpdfstring{$\phi$}{𝜙}-Assouad dimensions}\label{s:general-properties}
In this section, we prove properties of the $\phi$-Assouad dimensions which hold for general sets.
We recall the definitions introduced in the introduction: in particular, see \cref{ss:gen-Assouad} and \cref{ss:windows}.
\subsection{Moran sets, window bounds and topology}\label{ss:window-proofs}
In this section, we establish the bounds within windows as stated in \cref{it:fund-dim}.
In order to do so, we first need to recall the notion of a homogeneous Moran set in $\R^d$; we also study these sets in~\cref{ss:moran}.
Let $\mathcal{I}=\{0,1\}^d$, set $\mathcal{I}^*=\bigcup_{n=0}^\infty\mathcal{I}^n$, and denote the unique word of length $0$ by $\varnothing$.
Suppose the ratios $\bm{r}=(r_n)_{n=1}^\infty$ satisfy $0<r_n\leq 1/2$ for each $n\in\N$.
Then for $n\in\N$ and $\bm{i}\in\mathcal{I}$, we let $S^n_{\bm{i}}\colon\R^d\to\R^d$ denote the similarity
\begin{equation*}
    S^n_{\bm{i}}(x)\coloneqq r_n x+b^n_{\bm{i}}
\end{equation*}
where the $j$th entry of $b^n_{\bm{i}}\in\R^d$ is
\begin{equation*}
    (b^n_{\bm{i}})^{(j)} =
    \begin{cases}
      0 &\text{if } \bm{i}^{(j)}=0,\\*
      1-r_n &\text{if } \bm{i}^{(j)}=1.
    \end{cases}.
\end{equation*}
We extend this definition to finite words $\sigma=(\bm{i}_1,\ldots,\bm{i}_n)\in\mathcal{I}^n$ by $S_\sigma=S^1_{\bm{i}_1}\circ\cdots\circ S^n_{\bm{i}_n}$.
Finally, set
\begin{equation*}
    M_n=\bigcup_{\sigma\in\mathcal{I}^n}S_\sigma([0,1]^d)\qquad\text{and}\qquad M=M(\bm{r})\coloneqq\bigcap_{n=1}^\infty M_n.
\end{equation*}
We refer to the set $M$ as a \defn{homogeneous Moran set} (with contraction ratios $\bm{r}$).
When dealing with homogeneous Moran sets, it is convenient to work with the max norm.
Then $M_n$ consists of $2^{dn}$ hypercubes each with diameter $\rho_n\coloneqq r_1\cdots r_n$.
We will denote such a hypercube by $Q_n$, i.e.\ $Q_n=S_\sigma([0,1]^d)$ for some $\sigma\in\mathcal{I}^n$.
We now give a formula for the $\phi$-Assouad dimensions of homogeneous Moran sets.
\begin{proposition}\label{p:moran-dim-formula}
    Fix a homogeneous Moran set $M$ and dimension function $\phi$.
    Then
    \begin{align*}
        \dimAs{\phi} M &= \limsup_{n \to \infty} \frac{(m_n - n)d \log 2}{-\phi(\rho_n) \log \rho_n},\\
        \dimuAs{\phi} M &= \limsup_{n \to \infty} \sup_{m \geq m_n} \frac{(m - n)d \log 2}{\log (\rho_n / \rho_m)},
    \end{align*}
    where for all $n\in\N$,
    \begin{equation*}
        m_n \coloneqq \max\{ m \geq n : \rho_m > \rho_n^{1 + \phi(\rho_n)} \}.
    \end{equation*}
\end{proposition}
\begin{proof}
    We give a short proof for $\dimAs{\phi} M$; the proof for the upper $\phi$-Assouad dimension is similar.
    Note that for any level-$n$ hypercube $Q$,
    \begin{equation*}
        N_{\rho_n^{1 + \phi(\rho_n)}}(Q) = 2^{d(m_n-n+1)},
    \end{equation*}
    which upon taking logarithms and a limit supremum, and using the definition of the $\phi$-Assouad dimension, proves the lower bound.
    For the upper bound, if $B$ is any hypercube of side-length $R$, let $n$ be such that $\rho_n \geq R > \rho_{n+1}$.
    Then
    \begin{equation*}
        N_{\rho_n^{1 + \phi(\rho_n)}}(B \cap M) \leq 2^{d} N_{\rho_{m_n+1}}(Q) \leq 2^{d(m_n - n + 2)}.
    \end{equation*}
    Using the fact that $R > \rho_{n+1}$ and $\rho_n^{1 + \phi(\rho_n)} < \rho_{m_n}$, the upper bound follows.
\end{proof}

The next result will be used for the proof of \cref{it:fund-dim}~\cref{im:limit-diff}.
\begin{proposition}\label{p:continuity-converse}
    Fix $d \in \N$ and $\eps \in (0,1)$.
    If $\phi,\psi$ are any dimension functions such that $\liminf_{R \to 0} \frac{\phi(R)}{\psi(R)} < 1-\eps$, then there exists a homogeneous Moran set $M \subset \R^d$ such that $\dimAs\phi M = d$ and $\dimAs\psi M \leq d - d \eps$.
\end{proposition}
\begin{proof}
    By assumption, there is a sequence of scales $R_n$ converging to $0$ such that
    \begin{equation*}
        \frac{\phi(R_n)}{\psi(R_n)} < 1-\eps.
    \end{equation*}
    Without loss of generality we may assume
    \begin{equation*}
        R_{n+1} < \min\{R_n^{1+\psi(R_n)}/4,R_n^n\}
    \end{equation*}
    for all $n$.
    The idea is to construct a Moran set whose covering number increases as fast as possible between scales $R_n$ and $R_n^{1+\phi(R_n)}$, but decreases as fast as possible at all other scales.
    Indeed, let $r_1 = \rho_1 = R_1$, and let $k_1 = 1$.
    Assume that we have defined $r_1,r_2,\dotsc,r_{k_n}$ for some $n$, and that $\rho_{k_n} = R_n$.
    Let $k_{n+1}$ be the smallest integer such that
    \begin{equation*}
        \rho_{k_n}\cdot 2^{-(k_{n+1} - k_n - 1)} < R_n^{1+\phi(R_n)}.
    \end{equation*}
    Let $r_{k_n + 1} = \dotsb = r_{k_{n+1}-1} = 1/2$, and let $r_{k_{n+1}}$ be such that $\rho_{k_{n+1}} = R_{n+1}$.
    Note that $r_{k_{n+1}} < 1/2$ since $R_{n+1} < R_n^{1+\psi(R_n)}/4$.
    Let $M$ be the homogeneous Moran set obtained by this inductive process.

    Since $r_{j} = 1/2$ for all $j$ corresponding to scales between $R_n$ and $R_n^{1+\phi(n)}$,
    \begin{equation*}
        N_{R_n^{1+\phi(R_n)}}(M \cap [0,R_n]^d) \approx R_n^{-d\cdot\phi(R_n)}
    \end{equation*}
    with implicit constants independent of $n$, so $\dimAs\phi M = d$.
    On the other hand, since $\rho_{k_{n+1}} = R_{n+1} \leq R_n^{1+\psi(n)}$ by assumption, for any $x\in M$,
    \begin{equation*}
        N_{R_{n}^{1+\psi(n)}}(M \cap B(x, R_n^{1+\phi(n)})) \approx 1.
    \end{equation*}
    Therefore since $\psi$ is monotonic and $R_{n+1} < R_n^n$,
    \begin{align*}
        \dimAs\psi F &= \limsup_{R \to 0} \sup_{x \in \R^d} \frac{\log N_{R^{1+\psi(R)}} (M \cap B(x,R))}{-\psi(R) \log R} \\
                                &\leq \limsup_{n \to \infty} \frac{\log N_{R_n^{1+\psi(R_n)}} (M \cap [0,R_n]^d)}{-\psi(R_n) \log R_n}  \\
                                &\leq \limsup_{n \to \infty} \frac{- d \phi(R_n) \log R_n}{-\psi(R_n)\log R_n} \\
                                &\leq d-d\eps
    \end{align*}
    as claimed.
\end{proof}
Essentially by combing this result with \cite[Proposition~2.11]{zbl:1485.28006}, we obtain the following corollary.
We provide the details here for completeness since our notion of $\phi$-Assouad dimension is slightly different than the notion in their paper.
\begin{corollary}\label{c:phi-equivs}
    Let $\eps>0$ and suppose $\phi$ and $\psi$ are any functions from $(0,1)\to\R^+$ with
    \begin{equation*}
        1-\eps<\liminf_{R\to 0}\frac{\phi(R)}{\psi(R)}\leq\limsup_{R\to 0}\frac{\phi(R)}{\psi(R)}< 1+\eps.
    \end{equation*}
    Then for any n.b.d.\ space $F$ with doubling constant $M$,
    \begin{equation*}
        |\dimAs\psi F-\dimAs\phi F|\leq \eps(1+2\log_2 M+\eps).
    \end{equation*}
    In particular, suppose $\phi$ and $\psi$ are dimension functions.
    Then the following are equivalent:
    \begin{enumerate}[nl,r]
        \item\label{im:phi-eq-1} $\lim_{R\to 0}\phi(R)/\psi(R)=1$.
        \item\label{im:phi-eq-2} For all bounded $F \subset \R$, $\dimAs\psi F = \dimAs\phi F$.
        \item\label{im:phi-eq-3} For all n.b.d.\ spaces $F$, $\dimAs\psi F = \dimAs\phi F$.
    \end{enumerate}
\end{corollary}
\begin{proof}
    Let $\eps>0$ be arbitrary and let $R'$ be sufficiently small so that $|\phi(R)/\psi(R)-1|<\eps$ for all $R \in (0,R']$.
    Now let $C>0$ be such that for all $R\in(0,R')$ and $x\in F$,
    \begin{equation*}
        N_{R^{1+\phi(R)}}(F\cap B(x,R))\leq C\Bigl(\frac{R}{R^{1+\phi(R)}}\Bigr)^{\dimAs\phi F+\eps}.
    \end{equation*}
    Now for all $R\in(0,R')$,
    \begin{equation*}
        R^{1+\psi(R)}=R^{1+\phi(R)}R^{\psi(R)-\phi(R)}\geq R^{1+\phi(R)}R^{\psi(R)\eps}=R^{1+\phi(R)}2^{\psi(R)\eps\log_2R},
    \end{equation*}
    so each ball of radius $R^{1+\psi(R)}$ can be covered by at most $M \cdot R^{-\log_2 M\cdot \psi(R)\cdot\eps}$ balls of radius $R^{1+\phi(R)}$.
    Thus
    \begin{align*}
        N_{R^{1+\psi(R)}}(F\cap B(x,R))&\leq M\cdot R^{-\log_2 M\cdot \psi(R)\cdot\eps}\cdot N_{R^{1+\phi(R)}}(F\cap B(x,R))\\
                                                      &\leq M\cdot C\cdot R^{-\log_2 M\cdot \psi(R)\cdot\eps}\Bigl(\frac{R}{R^{1+\phi(R)}}\Bigr)^{\dimAs\phi F+\eps}\\
                                                      &\leq M\cdot C\cdot \Bigl(\frac{R}{R^{1+\psi(R)}}\Bigr)^{(1+\eps)(\dimAs\phi F +\eps)+\eps\log_2 M}.
    \end{align*}
    Therefore
    \begin{equation*}
        \dimAs\psi F\leq(1+\eps)(\dimAs\phi F+\eps)+\eps\log_2 M.
    \end{equation*}
    and since $\dimAs\phi F\leq\dimA F\leq\log_2 M$,
    \begin{equation*}
        \dimAs\psi F-\dimAs\phi F\leq\eps(1+2\log_2 M+\eps).
    \end{equation*}
    The reverse inequality follows by the same argument, as required.

    This also proves that \cref{im:phi-eq-1} implies \cref{im:phi-eq-3}.
    To see the remaining equivalences, the implication \cref{im:phi-eq-3} implies \cref{im:phi-eq-2} is immediate and that \cref{im:phi-eq-2} implies \cref{im:phi-eq-1} follows from \cref{p:continuity-converse}.
\end{proof}
Next, we observe the following analogue, for arbitrary dimension functions converging to $0$, of the usual bounds for the Assouad spectrum given in \cite[Proposition~3.4]{zbl:1390.28019}.
This result is a mild specialization and refinement of \cite[Proposition~2.15]{zbl:1485.28006}, and follows by a similar strategy.
\begin{proposition}\label{p:window-bounds}
    Let $\phi$ be a dimension function and suppose $\lim_{R\to 0}\phi(R)=0$.
    Let $F$ be a n.b.d.\ space and let $\varphi(\gamma)\coloneqq \dimAs{\phi_\gamma} F$.
    Then for all $0<\alpha<\beta<\infty$,
    \begin{equation}\label{e:ab-bounds}
        0\leq\frac{1}{\alpha}\varphi(\alpha)-\frac{1}{\beta}\varphi(\beta)\leq\frac{\beta-\alpha}{\alpha\beta}\varphi\left(\frac{\alpha\beta}{\beta-\alpha}\right).
    \end{equation}
    In particular, $\varphi$ is a continuous function of $\gamma$.
\end{proposition}
\begin{proof}
    Let $0<\alpha<\beta<\infty$ and $\eps>0$ be arbitrary.
    Let $R\in(0,1)$ and $x\in F$.
    Since
    \begin{equation*}
        B\bigl(x,R^{\frac{\beta}{\beta+\phi(R)}}\bigr)\subseteq B\bigl(x,R^{\frac{\alpha}{\alpha+\phi(R)}}\bigr),
    \end{equation*}
    it follows that for all $R$ sufficiently small
    \begin{align*}
        \sup_{x\in F}N_R\bigl(F\cap B\bigl(x,R^{\frac{\alpha}{\alpha+\phi(R)}}\bigr)\bigr)&\geq \sup_{x\in F}N_R\bigl(F\cap B\bigl(x,R^{\frac{\beta}{\beta+\phi(R)}}\bigr)\bigr)\\
                                                                                                         &\geq \Bigl(\frac{R^{\frac{\beta}{\beta+\phi(R)}}}{R}\Bigr)^{\varphi(\beta)-\eps}\\
                                                                                                         &= \Bigl(\frac{R^{\frac{\alpha}{\alpha+\phi(R)}}}{R}\Bigr)^{\frac{\alpha+\phi(R)}{\beta+\phi(R)}\cdot(\varphi(\beta)-\eps)}.
    \end{align*}
    Since $\lim_{R\to 0}\phi(R)=0$ and $\eps>0$ was arbitrary, the first inequality in \cref{e:ab-bounds} follows.

    We obtain the second inequality by covering balls of radius $R^{\frac{\alpha}{\alpha+\phi(R)}}$ with balls of radius $R^{\frac{\beta}{\beta+\phi(R)}}$, to give
    \begin{align*}
        \sup_{x\in F}N_{R}(F \cap B(x,R^{\frac{\alpha}{\alpha+\phi(R)}}))\leq{}&\sup_{x\in F}N_{R^{\frac{\beta}{\beta+\phi(R)}}}(F \cap B(x,R^{\frac{\alpha}{\alpha+\phi(R)}}))\\*
                                                                              &\cdot\sup_{x\in F} N_R(F \cap B(x,R^{\frac{\beta}{\beta+\phi(R)}})).
    \end{align*}
    Next, observe that
    \begin{equation*}
      \frac{R^{\frac{\alpha}{\alpha+\phi(R)}}}{R^{\frac{\beta}{\beta+\phi(R)}}}=\frac{R}{R^{1+\eta(R)}},
    \end{equation*}
    where $\eta(R)=\frac{\beta}{\beta+\phi(R)}-\frac{\alpha}{\alpha+\phi(R)}$, and
    \begin{align*}
        \lim_{R\to 0}\frac{\phi(R)}{\eta(R)}=\lim_{R\to 0}\frac{(\alpha+\phi(R))(\beta+\phi(R))}{\beta-\alpha}=\frac{\alpha\beta}{\beta-\alpha}.
    \end{align*}
    In particular, by \cref{c:phi-equivs}, $\dimAs\eta F=\varphi\bigl(\frac{\alpha\beta}{\beta-\alpha}\bigr)$.
    Fix $\eps>0$.
    Then for all $R$ sufficiently small,
    \begin{align*}
        \sup_{x\in F}N_R(F \cap B\bigl(x,R^{\frac{\alpha}{\alpha+\phi(R)}}\bigr))&\leq\Bigl(\frac{R^{\frac{\alpha}{\alpha+\phi(R)}}}{R^{\frac{\beta}{\beta+\phi(R)}}}\Bigr)^{\varphi\left(\frac{\alpha\beta}{\beta-\alpha}\right)+\eps}\cdot\Bigl(\frac{R^{\frac{\beta}{\beta+\phi(R)}}}{R}\Bigr)^{\varphi(\beta)+\eps}\\
                                                                                               &\leq\Bigl(\frac{R^{\frac{\alpha}{\alpha+\phi(R)}}}{R}\Bigr)^{\frac{\alpha+\phi(R)}{-\phi(R)}\bigl(-\eta(R)\left(\varphi\left(\frac{\alpha\beta}{\beta-\alpha}\right)+\eps\right)+\frac{-\phi(R)}{\beta+\phi(R)}(\varphi(\beta)+\eps)\bigr)}\\
                                                                                               &=\left(\frac{R^{\frac{\alpha}{\alpha+\phi(R)}}}{R}\right)^{\frac{\beta-\alpha}{\beta+\phi(R)}\left(\varphi\left(\frac{\alpha\beta}{\beta-\alpha}\right)+\eps\right)+\frac{\alpha+\phi(R)}{\beta+\phi(R)}(\varphi(\beta)+\eps)}
    \end{align*}
    and taking the limit as $R$ goes to zero yields the desired bound.

    It is immediate that $\varphi$ is a continuous function since $\varphi(\gamma)\leq\dimA F<\infty$.
\end{proof}
\begin{remark}
    The assumption that $\lim_{R\to 0}\phi(R)=0$ in \cref{p:window-bounds} is precisely saying that the dimension function $\phi$ is not equivalent to the Assouad spectrum at some $\theta\in(0,1)$ by \cref{it:fund-dim}~\cref{im:limit-eq} (recalling \cref{ex:spectrum-equiv}).
    The same proof works if instead $\lim_{R\to 0}\phi(R)=\theta$ (though the resulting formula is slightly different), so \cref{it:fund-dim}~\cref{im:cont} does indeed hold.
    These details are proved in the same way as the usual bounds for the Assouad spectrum given in \cite[Proposition~3.4]{zbl:1407.28002}.
\end{remark}
Our next result in particular implies the general bounds between disjoint windows given in \cref{it:fund-dim}~\cref{im:window-compare}.
\begin{proposition}\label{p:inter-window-compare}
    Let $F$ be a n.b.d.\ space and let $\psi$ be a dimension function.
    Then for all $\eps>0$ there exists $\eta>0$ (depending on $\dimA F$, $\psi$ and $\eps$) such that if $\phi$ is any dimension function with $\limsup_{R\to 0}\frac{\phi(R)}{\psi(R)}< \eta$ then
    \begin{equation*}
        \dimuAs{\psi} F \leq \dimAs{\phi} F+\eps.
    \end{equation*}
\end{proposition}
\begin{proof}
    The idea is to cover the intersection of $F$ with a ball with smaller balls corresponding to the scale given by $\phi$, and then similarly cover the intersection of $F$ with each of those balls, continuing until we reach approximately the desired scale.
    When $\phi(R)/\psi(R)$ is small, the error resulting from not hitting the exact scale will be negligible.

    Fix $\eps>0$.
    Since $\psi$ is bounded, using L'Hôpital's rule there exists $\eta>0$ such that
    \begin{equation*}
        \eta < \frac{(1-\frac{1}{1+\psi(R)})\eps}{\psi(R)(s+\eps)}
    \end{equation*}
    for all $R \in (0,1)$ and $s \coloneqq \dimAs{\phi} F\in[0,\dimA F]$.
    By assumption, there exists $\overline{R}>0$ such that for all $0<R\leq\overline{R}$,
    \begin{enumerate}[nl]
        \item $\phi(R)/\psi(R) < \eta$, and
        \item\label{im:window-sup-cover} if $B$ is any ball of radius $R$ intersecting $F$ then $N_{R^{1+\phi(R)}}(F \cap B)\leq R^{-\phi(R)(s+\eps)}$.
    \end{enumerate}
    Now fix $0<R\leq \overline{R}$ and $r\leq R^{1+\psi(R)}$.
    Define the strictly decreasing sequence $(R_n)_{n \in \N}$ by $R_0 = R$ and $R_n = R_{n-1}^{1+\phi(R_{n-1})}$ for $n \in \N$.
    Let $n(R)$ be the smallest natural number such that $R_{n(R)}\leq r$.
    Applying \cref{im:window-sup-cover} inductively gives
    \begin{align*}
        N_r(F \cap B(x,R)) &\leq N_{R_{n(R)}}(F \cap B(x,R)) \\
        &\leq \prod_{i=1}^{n(R)} \left( \frac{R_{i-1}}{R_i} \right)^{s+\eps} \\
                           &\leq \left(\frac{R}{r}\right)^{s+\eps} \left(\frac{R_{n(R)-1}}{R_{n(R)}}\right)^{s+\eps} \\
                           &\leq \left(\frac{R}{r}\right)^{s+\eps} r^{-\phi(R)(s+\eps)} \\
                           &\leq \left(\frac{R}{r}\right)^{s+2\eps},
    \end{align*}
    where the last line follows by the choice of $\eta$.
\end{proof}
We obtain the following result as a direct application.
\begin{corollary}\label{c:dimub-lower}
    Let $\phi$ be any dimension function and $F$ a n.b.d.\ space.
    Then $\dimuB F\leq\dimAs\phi F$.
\end{corollary}
Note that there exists $F$ for which $\dimuB F < \dimAs\phi F$ for all dimension functions $\phi$; see \cref{r:no-upper-box-dim-func}.
Next, we provide the proof of the characterization of the topology on the space of dimension functions.
Recall that $f_F([\psi]) = \dimAs\psi F$.
\begin{proofref}{ic:topology}
    Let $B \subset \R$ be open and $F \subset \R^d$ be bounded.
    Suppose $\psi$ is such that $\dimAs\psi F \in B$.
    Then by \cref{it:fund-dim}~\cref{im:cont}, for all $\eps > 0$ sufficiently small, $\dimAs{(1+\eps)\psi} F, \dimAs{(1-\eps)\psi} F \in B$.
    Therefore $N_{\psi,\eps} \subset f_F^{-1}(B)$.
    Thus $f_F^{-1}(B)$ is open with respect to the topology generated by the $N_{\psi,\eps}$.

    Now fix a dimension function $\psi$ and $\eps \in (0,1)$.
    Suppose $\phi$ is such that
    \begin{equation*}
        1-\eps < \liminf_{R \to 0} \frac{\phi(R)}{\psi(R)} \leq \limsup_{R \to 0} \frac{\phi(R)}{\psi(R)} < 1+\eps.
    \end{equation*}
    By \cref{it:fund-dim} \cref{im:limit-diff}, there exist bounded $M_1,M_2 \subset \R$ and $\eps'>0$ such that
    \begin{equation*}
        \dimAs{(1-\eps)\psi} M_1 = \dimAs{\phi} M_2 = d
    \end{equation*}
    and
    \begin{equation*}
        \max\{\dimAs{\phi} M_1 , \dimAs{(1+\eps)\psi} M_2\} < d-\eps'.
    \end{equation*}
    Then
    \begin{equation*}
        \phi \in f_{M_1}^{-1}((-1,d)) \cap f_{M_2}^{-1}((d-\eps',d+\eps')) \subseteq N_{\psi,\eps}.
    \end{equation*}
    Therefore $N_{\psi,\eps} \in \mathcal{T}$.
\end{proofref}
To conclude this section, we establish a general result on the existence of maximal dimension functions.
This result implies that the family of dimension functions is quite rich, and will also be useful in the proof of \cref{t:Assouad-recover}.
\begin{proposition}\label{p:dim-maximal}
    Let $\phi\colon(0,1)\to(0,1)$ be a continuous function such that
    \begin{equation*}
        \lim_{R\to 0}\phi(R)\log(1/R)=\infty.
    \end{equation*}
    Then the set of dimension functions
    \begin{equation*}
        \{\psi\in\mathcal{W}:\psi\leq\phi\}
    \end{equation*}
    has a unique maximal element with respect to pointwise comparison.
\end{proposition}
\begin{proof}
    First, let
    \begin{equation*}
        \psi_0(R)=\frac{\inf_{0<r\leq R}\phi(r)\log(1/r)}{\log(1/R)}.
    \end{equation*}
    By definition, $\psi_0(R)\leq\phi(R)$ and $\psi_0(R)\log(1/R)$ increases to infinity as $R$ decreases to zero, and moreover $\psi_0$ is the unique maximal function with these properties.

    Next, define
    \begin{equation*}
        \psi(R)=\inf_{r\in(R,1)}\psi_0(R).
    \end{equation*}
    Of course, $\psi(R)\leq\psi_0(R)$ and $\phi$ is monotonically decreasing, and since $\phi$ is continuous, $\psi(R)>0$ for all $R\in(0,1)$.
    Moreover, $\psi$ is the unique maximal function $\psi\leq\psi_0$ with these properties.

    It remains to prove that $\psi(R)\log(1/R)$ increases to infinity as $R$ decreases to zero.
    First, for any $0<R<1$, either $\psi(R)=\psi_0(R)$ (and we set $r_a(R)=R=r_b(R)$) or there are $r_a(R)<R<r_b(R)$ such that $\psi(r_a(R))=\psi_0(r_a(R))=\psi(r_b(R))=\psi_0(r_b(R))$ so $\psi$ is constant on $[r_a(R),r_b(R)]$.
    Suppose $0<R_0<R_1<1$, and we may assume that $r_b(R_0) < r_a(R_1)$.
    Then using the properties of $\psi_0$,
    \begin{align*}
        \psi(R_0)\log(1/R_0)=\psi_0(r_b(R_0))\log(1/R_0)&\geq\psi_0(r_b(R_0))\log(1/r_b(R_0))\\
                                             &\geq\psi_0(r_a(R_1))\log(1/r_a(R_1))\\
                                             &\geq\psi(r_a(R_1))\log(1/R_1)\\
                                             &=\psi(R_1)\log(1/R_1),
    \end{align*}
    as required.
\end{proof}
\subsection{Recovering the upper \texorpdfstring{$\phi$}{𝜙}-Assouad dimension}\label{ss:upper-recover}
Next, we obtain a generalization of \cite[Theorem~2.1]{zbl:1410.28008}, in the case where $F$ is a bounded set.
This result is stated in \cref{it:upper-recover} and gives a formula for the upper $\phi$-Assouad dimension in terms of the $\psi$-Assouad dimensions of functions in the same rate window as $\phi$.
\begin{proofref}{it:upper-recover}
    The result is trivial if $s\coloneqq\dimuAs\phi F=0$.
    Otherwise, if $\lim_{R\to 0}\phi(R)>0$, then $\dimAs\phi F$ is the Assouad spectrum of $F$ at $\theta$ where ${\frac{1}{\theta}-1}=\lim_{R\to 0}\phi(R)$ by \cref{c:phi-equivs}.
    In particular, the result follows by \cite[Theorem~2.1]{zbl:1410.28008}.
    Thus we may assume $\lim_{R\to 0}\phi(R)=0$.
    Moreover, since $\phi_\alpha\geq\phi$ for all $\alpha\in(0,1)$, it always holds that $\dimuAs\phi F\geq\dimAs{\phi_\alpha} F$ for all $\alpha\in(0,1)$.
    It suffices to show the converse inequality.

    Let $(\eps_n)_{n=1}^\infty\subset(0,s)$ converge monotonically to $0$.
    By definition, we can find a sequence $(x_n,r_n,R_n)_{n=1}^\infty$ such that $x_n\in F$, $0<r_n\leq R_n^{1+\phi(R_n)}<R_n\leq 1$, $\lim_{n\to\infty}r_n/R_n=0$, and
    \begin{equation}\label{e:eps-bd}
        \Bigl(\frac{R_n}{r_n}\Bigr)^{s-\eps_n}\leq N_{r_n}(F \cap B(x_n,R_n)).
    \end{equation}
    For each $n$, let $\alpha_n$ be such that $R_n^{1+\phi_{\alpha_n}(R_n)}=r_n$, and observe that $\alpha_n\in(0,1]$.
    If $\dimuB F\geq s$ then $\dimuAs\phi F=\dimAs{\phi_\alpha}F$ for all $\alpha\in(0,1)$ by \cref{c:dimub-lower} and we are done, so we may assume that $R_n \to 0$ monotonically as $n \to \infty$.
    Let $b_n=\phi(R_n)/\alpha_n$.
    Passing to a subsequence if necessary, we may assume that either $(b_n)_{n=1}^\infty$ diverges to infinity, or $(\alpha_n)_{n=1}^\infty$ and $(b_n)_{n=1}^\infty$ respectively converge monotonically to $\alpha_0\in[0,1]$ and $b_0\in[0,\infty)$.

    If $(b_n)_{n=1}^\infty$ diverges to infinity, let $\psi$ denote the constant function $1$, and if $b_0\in(0,\infty)$, let $\psi$ denote the constant function $b_0/2$.
    In either case, \cref{e:eps-bd} and \cref{p:inter-window-compare} imply that $s\leq\dimuAs\psi F\leq\dimAs\phi F$, so we are done by \cref{p:window-bounds}.

    Otherwise, $b_n$ decreases to $0$.
    If $\alpha_0=0$, since $b_n\geq\phi(R_n)$, we may choose $\psi$ to be a monotonically decreasing function such that $\psi(R_n)=b_n$ and $\psi(R)\geq\phi(R)/\alpha_n$ for all $R\in[R_{n+1},R_n]$.
    It is clear that $\psi$ is a dimension function since the $\alpha_n$ are monotonically decreasing.
    Moreover, since $\psi(R_n)=\phi(R_n)/\alpha_n$, by \cref{e:eps-bd}, $s\leq\dimAs\psi F\leq\dimuAs\psi F$.
    Finally, since $\alpha_0=0$, $\lim_{R\to0}\phi(R)/\psi(R)=0$, so that $s\leq\dimuAs\psi F\leq\dimAs\phi F$ by \cref{p:inter-window-compare}.

    In the final case, $\alpha_0>0$, and since $b_n$ decreases to $0$, we may choose $\psi$ to be a function such that $\psi(R)\leq R$, $\psi(R_n)=b_n$, $\lim_{R\to 0}\phi(R)/\psi(R)=\alpha_0$.
    Again, \cref{e:eps-bd} implies that $s\leq\dimAs\psi F$, and it follows by \cref{c:phi-equivs} that $\dimAs\psi F=\dimAs{\phi_{\alpha_0}} F$.
    Thus $\dimAs{\phi_{\alpha_0}}F=s$.
\end{proofref}

\subsection{Recovering the interpolation}\label{ss:interpolation-recover}
In this section, we prove that the $\phi$-Assouad dimensions \emph{recover the interpolation}.
We first show that the Assouad dimension is attained as the $\phi$-Assouad dimension for some dimension function $\phi$.
Moreover, $\phi$ can be chosen to be arbitrarily small, implying that in the definition of the Assouad dimension it suffices to consider scales $r$ and $R$ which are very close together.
\begin{theorem}\label{t:Assouad-recover}
    Let $g\colon(0,1)\to(0,1)$ be a continuous function such that
    \begin{equation*}
        \lim_{R\to 0}g(R)/R=0,
    \end{equation*}
    and let $F$ be any n.b.d.\ space.
    Then there is a dimension function $\psi$ with $R^{1+\psi(R)}\geq g(R)$ for all $R\in(0,1)$ such that
    \begin{equation*}
        \dimA F=\dimuAs\psi F=\dimAs\psi F=\limsup_{R\to 0}\frac{\log \sup_{x\in F}N_{R^{1+\psi(R)}}\bigl(F \cap B(x,R)\bigr)}{-\psi(R)\log R}.
    \end{equation*}
\end{theorem}
\begin{proof}
    Write $\phi(R)=\frac{\log g(R)}{\log R}-1$ and note that $\phi(R)\log(1/R)$ increases to infinity as $R$ decreases to zero by the assumption on $g$.
    Thus applying \cref{p:dim-maximal}, there is a unique maximal dimension function $\psi_0(R)\leq\phi(R)$.
    The same proof as \cref{p:inter-window-compare} gives that for any function $h\colon(0,1)\to(0,1)$ satisfying $h(R)\leq R$ and $h(R)/R \to 0$ as $R \to 0$,
    \begin{align*}
        \dimA F=\inf\Bigl\{s:(\exists C>0)&(\forall 0<h(R)\leq r\leq R<1)\\*
                                                  &\sup_{x\in F}N_r(F\cap B(x,R))\leq C \Bigl(\frac{R}{r}\Bigr)^s\Bigr\};
    \end{align*}
    in particular, this holds for $h(R)=R^{1+\psi_0(R)}$.
    Since $F$ is doubling, get a sequence $(R_n,r_n,x_n)_{n=1}^\infty$ such that $x_n\in F$, $R_n$ and $r_n/R_n$ decrease monotonically to 0, $r_n\geq h(R_n)$, and
    \begin{equation*}
        \dimA F=\lim_{n\to\infty}\frac{\log N_{r_n}(F\cap B(x_n,R_n))}{\log(R_n/r_n)}.
    \end{equation*}
    For each $n\in\N$, let $\theta_n$ be such that $r_n=R_n^{1+\theta_n}$.
    Note that $\theta_n\leq\psi_0(R_n)$ by the assumption on $r_n$.
    Since $r_n/R_n$ decreases monotonically to $0$, $\theta_n\log(1/R_n)$ diverges monotonically to infinity.
    Since $\psi_0(R)$ decreases as $R \to 0$, passing to a subsequence, $\theta_n$ converges to some $\theta\in[0,\infty)$.
    If $\theta>0$, the function $\psi_1$ defined to take the constant value $\theta_n$ on each interval $(R_{n+1},R_n]$ has $\dimAs{\psi_1}F=\dimA F$.
    Of course, $\psi_1$ need not be a dimension function, but if $\psi$ denotes the constant function $\theta$, by \cref{c:phi-equivs}, $\dimAs\psi F=\dimAs{\psi_1}F=\dimA F$.
    Since $\theta\leq\lim_{R\to 0}\phi(R)$, necessarily $\psi\leq\psi_0$, so $\psi$ satisfies the required properties.

    Otherwise, again passing to a subsequence, we may assume that $\theta_n$ decreases strictly to $0$.
    It suffices to choose a dimension function $\psi$ such that $\psi(R)\leq \psi_0(R)$ and $R_n^{1+\psi(R_n)}=r_n$.
    Indeed, assuming we have found such a function $\psi$, recalling the formula in \cref{e:phi-dimA-lim},
    \begin{align*}
        \dimA F &\geq\dimuAs\psi F\\
                &\geq\dimAs\psi F\\
                &=\limsup_{R\to 0}\frac{\log \sup_{x\in F}N_{R^{1+\psi(R)}}\bigl(F \cap B(x,R)\bigr)}{-\psi(R)\log R}\\
                &\geq\limsup_{n\to\infty}\frac{\log N_{R_n^{1+\psi(R_n)}}\bigl(F \cap B(x_n,R_n)\bigr)}{-\psi(R_n)\log R_n}\\
                &=\lim_{n\to\infty}\frac{\log N_{r_n}(F\cap B(x_n,R_n))}{\log(R_n/r_n)}\\
                &= \dimA F
    \end{align*}
    so that all the inequalities are equalities, as claimed.

    We now inductively define such a function $\psi$ on $(0,1)$ as follows.
    Define $\psi$ on the interval $[R_1,1)$ to be the constant function $\theta_1$ and note that $\psi$ satisfies the conditions of being a dimension function on $[R_1,1)$.
    Suppose $\psi$ is defined on the interval $[R_n,1)$ for some $n\in\N$.
    Let $R_{n+1}' < R_n$ be such that
    \begin{equation*}
        \theta_n\frac{\log R_n}{\log R_{n+1}'}=\theta_{n+1}
    \end{equation*}
    and define $\psi$ to be the function $\theta_n\frac{\log R_n}{\log R}$ on $[R_{n+1}',R_n]$ and the constant function on $[R_{n+1},R_{n+1}']$.
    Note that $R_{n+1}'$ is chosen precisely so that $\psi$ is continuous at $R_{n+1}'$.
    Since $\theta_n\log(1/R_n)$ is increasing in $n$, $R_{n+1}\leq R_{n+1}'$.
    A direct check gives that $\psi$ satisfies the conditions of being a dimension function on $[R_{n+1},R_n]$.

    Finally, by construction, $\psi$ is the smallest possible choice of dimension function satisfying $\psi(R_n)=\theta_n$, and since $\theta_n\leq\psi_0(R_n)$, it follows that $\psi\leq\psi_0$.
\end{proof}
We can now complete the proof of the remaining cases in \cref{it:interpolation-recover}.
\begin{proofref}{it:interpolation-recover}
    If $\dimuB F < \alpha \leq \dimAs{\theta_0} F$ for some $\theta_0\in(0,1)$, then the constant function $1/\theta - 1$ where $\theta \coloneqq \inf\{ \theta' \in (0,1) : \dimAs{\theta'} F = \alpha \}$ satisfies the desired properties.
    We may thus assume that $\alpha>\dimuAs\theta F$ for all $\theta\in(0,1)$.
    Moreover, the case $\alpha=\dimA F$ is covered in \cref{t:Assouad-recover}, so we may assume that $\alpha<\dimA F$.

    The idea in the construction is for $\phi$ to remain constant for a long time until the dimension looks small, before decreasing at the fastest possible rate (while still satisfying the constraints of being a dimension function) until a carefully chosen scale.
    This process is then repeated inductively.
    First, let $\Omega=\{(R,r):0<r\leq R\leq 1\}$ and let $\omega\colon\Omega\to\R_{\geq 0}$ be given by
    \begin{equation*}
        \omega(R,r)=\frac{\sup_{x\in F}\log N_r\bigl(F \cap B(x,R)\bigr)}{\log (R/r)}.
    \end{equation*}
    For $(r,R)\in\Omega$, define
    \begin{align*}
        \mathcal{R}(R,r)&\coloneqq\bigl\{(x,y)\in\Omega:x\leq R,y\leq x^{1/\theta}\text{ where }r=R^{1/\theta}\bigr\}\\*
        \mathcal{A}(R,r)&\coloneqq\bigl\{(x,y)\in\Omega: x \leq R, x/y\geq R/r\bigr\}.
    \end{align*}
    Of course, $\mathcal{R}(R,r)\subseteq \mathcal{A}(R,r)$.
    We use the function $\omega$ and the regions $\mathcal{R}$ and $\mathcal{A}$ to define the decreasing sequences $(R_n)_{n \in \N}$, $(R_n')_{n \in \N}$ and the increasing sequence $(\theta_n)_{n \in \N}$ inductively as follows.

    Set $R_1 = 1$ and $\theta_1 = 1/2$, and assume we have defined $R_1,\ldots, R_n$, $R_1',\ldots,R_{n-1}'$ and $\theta_1,\ldots,\theta_n$ for some $n \in \N$.
    Since $\alpha > \dimuAs{\theta_n} F$, we can define the positive number
    \begin{equation*}
        R_n'\coloneqq\frac{1}{2}\sup\bigl\{R<R_n:\sup_{(x,y)\in\mathcal{R}(R,R^{1/\theta_n})}\omega(x,y)\leq \alpha\bigr\}.
    \end{equation*}
    Let
    \begin{align*}
        \mathcal{B}_n(R)\coloneqq{}&\bigcup_{R\leq\rho\leq R_n'}\mathcal{R}\bigl(\rho,\rho\cdot(R_n')^{\frac{1}{\theta_n}-1}\bigr)\\
        ={}&\mathcal{R}\bigl(R,R\cdot(R_n')^{\frac{1}{\theta_n}-1}\bigr)\cup\left(\mathcal{A}(R_n',(R_n')^{1/\theta_n})\cap\{(x,y):x\geq R\}\right).
    \end{align*}
    Since $\alpha < \dimA F$, $\omega(x,y)>\alpha$ for some $(x,y)\in\mathcal{A}(R_n',(R_n')^{1/\theta_n})\setminus\mathcal{R}(R_n',(R_n')^{1/\theta_n})$.
    Note that $(x,y) \in \mathcal{R}\bigl(x,x \cdot(R_n')^{\frac{1}{\theta_n}-1}\bigr) \subseteq \mathcal{B}_n(R)$ for all $R \leq x$.
    Therefore the number $R_{n+1}$ defined in the following way must be at least $x$, so in particular is positive:
    \begin{equation*}
        R_{n+1}\coloneqq\inf\Bigl\{R\leq R_n':\sup_{(x,y)\in\mathcal{B}_n(R)}\omega(x,y)\leq \alpha\Bigr\}.
    \end{equation*}
    \begin{figure}[t]
        \centering
        \begin{tikzpicture}[>=stealth,scale=10]
    \draw[->] (-0.01,0) -- (1.1,0);
    \draw[->] (0,-0.01) node[below]{$0$} -- (0,0.75);

    \matrix [draw, fill=white, below right,scale=0.8,column sep=0.4em, inner sep=0.4em, align=left] at (current bounding box.north west) {
        \draw[fill=black, fill opacity=0.15] (0,-0.32) rectangle (0.8,0); & \node[inner sep=0pt] {$\mathcal{B}_n(R_{n+1})$}; \\
        \fill[fill=black, opacity=0.07] (0,-0.40) -- (0.3, -0.4) -- (0.5,-0.08) -- (0, -0.08) -- cycle;
        \draw[fill=black, fill opacity=0.15] (0.8,-0.40) -- (0.3, -0.4) -- (0.5,-0.08) -- (0.8, -0.08) -- cycle;
        \draw[dashed] (0.5,-0.08) -- (0, -0.08) -- (0,-0.40) -- (0.3, -0.4);
        \draw (0.5,-0.08) -- (0.3, -0.4);
                                                                    & \node[inner sep=0pt] {$\mathcal{B}_n(\widetilde{R})\text{ for some }\widetilde{R}<R_{n+1}$}; \\
    };

    \coordinate (infl) at (0.8, 0.5724334);
    \coordinate (top) at (1,0.71554175);
    \coordinate (bottom) at (0.4,0.10119289);

    \fill[opacity=0.09, domain=0:0.8, smooth, variable=\t, samples=80, thick] plot ({\t}, {(\t)^(2.5)}) -- (top) -- (1,0) -- cycle;
    \fill[opacity=0.05,domain=0:0.7, smooth, variable=\t, samples=80, dashed] plot ({\t}, {(\t)^(1.93843242)}) -- (top) -- (1,0) -- cycle;

    \coordinate (omega) at (0.4,0.13523928);
    \node[circle, fill=black, inner sep=0.9pt, draw=black] at (omega) {};
    \node[left] at (omega) {$\omega(x,y)>\alpha$};

    \draw[dotted] (0,0) -- (infl);

    \draw[dotted] (0.8,-0.01) node[below]{$R_{n+1}$} -- (infl);
    \draw[dotted] (1,-0.01) node[below]{$R_{n}'$} -- (top);
    \draw[dotted] (0.2,-0.01) node[below]{$R_{n+1}'$} -- (0.2,0.01788854);

    \draw[domain=0:0.7, smooth, variable=\t, samples=80, dashed] plot ({\t}, {(\t)^(1.93843242)}) -- (top);
    \draw[domain=0:0.85, smooth, variable=\t, samples=30, thick] plot ({\t}, {(\t)^(2.5)}) node[above]{$y=x^{1/\theta_{n+1}}$};
\end{tikzpicture}
        \caption{Depiction of the choice of $R_{n+1}$ in the proof of \cref{it:interpolation-recover}.}
        \label{f:interpolation-image}
    \end{figure}
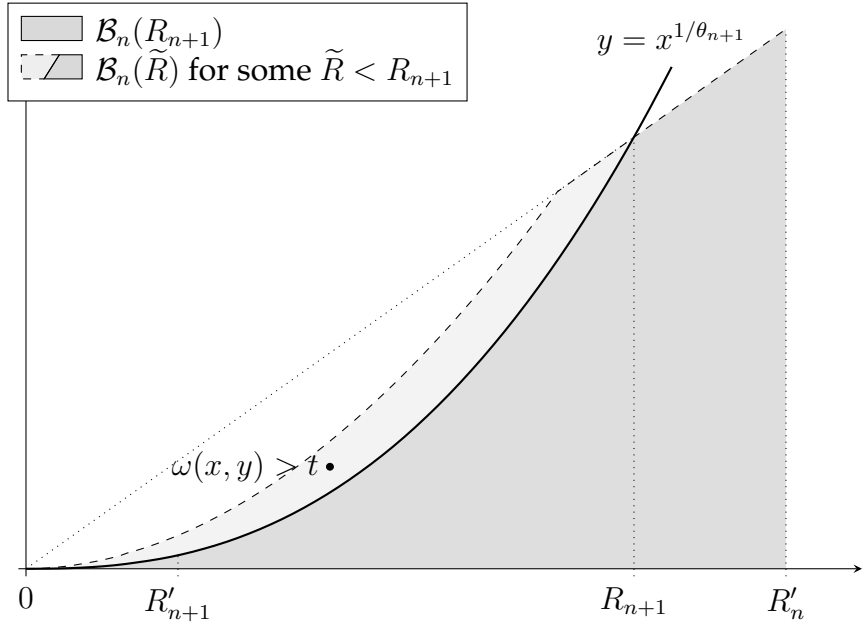
    The choice of $R_{n+1}$ is depicted in \cref{f:interpolation-image}.
    Define $\theta_{n+1}$ so that $R_{n+1}^{1/\theta_{n+1}}=R_{n+1}\cdot (R_n')^{\frac{1}{\theta_n}-1}$, or equivalently
    \begin{equation*}
        \theta_{n+1} \coloneqq \left( 1 + \left(\frac{1}{\theta_n} - 1 \right) \frac{\log R_n'}{\log R_{n+1}} \right)^{-1}.
    \end{equation*}
    If $R_n' = R_{n+1}'$ for some $n \in \N$, we can remove both $R_n'$ and $R_{n+1}'$ from the sequences, so we may assume that $R_1 > R_1' > R_2 > R_2' > \cdots$, and that $(\theta_n)_{n \in \N}$ is a strictly increasing sequence converging to $1$.
    We can now define $\phi$ by $\phi(x) = 1/\theta_n - 1$ for $x \in [R_n',R_n]$, and $\phi(x) \log x$ constant on $[R_{n+1},R_n']$, for all $n \in \N$.
    By construction, $\phi$ is a dimension function.

    It remains to prove dimension estimates.
    We begin with the upper bound.
    For $n \in \N$, if $R \in (R_{n+1},R_n')$ then by the definition of $R_{n+1}$, if $x \in F$, $r \leq R^{1+\phi(R)}$ then $N_r(B(x,R) \cap F) \leq (R/r)^\alpha$.
    Now suppose $R \in [R_n',R_n]$ for some $n > 1$, and let $x \in F$ and $r \leq R^{1+\phi(R)}$.
    Then
    \begin{equation*}
        N_r\bigl(B(x,R) \cap F\bigr) \leq N_r\bigl(B(x,2R) \cap F\bigr) \leq \Bigl(\frac{2R}{r}\Bigr)^\alpha,
    \end{equation*}
    where the last inequality follows from the definition of $R_n'$.
    It follows that $\dimuAs\phi F \leq \alpha$.
    
    We now prove that $\dimAs\phi F\geq \alpha$.
    For each $n > 1$, fix a positive number $\delta_n$ small enough that $\delta_n \leq (\theta_{n+1} - \theta_n)/n$ and
    \begin{equation*}
        1 - \frac{1}{\theta_n + \delta_n} \leq \left(1-\frac{1}{\theta_{n}}\right)\left(1-\frac{1}{n}\right).
    \end{equation*}
    By the definition of $R_{n}$, there exist $R \leq R_{n}$, $r \leq R^{\frac{1}{\theta_{n} + \delta_n}}$ and $x \in F$ such that $N_r (F \cap B(x,R)) > (R/r)^{\alpha}$.
    Then $R > R_{n+1}$ and $r \geq R^{1+\phi(R)}$.
    Now,
    \begin{align*}
        N_{R^{1+\phi(R)}}\bigl(B(x,R) \cap F\bigr)&\geq N_r\bigl(B(x,R) \cap F\bigr)\\
                                                  &> \Bigl(\frac{R}{r}\Bigr)^{\alpha}\\
                                                  &\geq R^{\left(1-\frac{1}{\theta_{n} + \delta_n}\right)\alpha} \\
                                                  &\geq R^{(1-1/\theta_{n})(1-1/n)\alpha}.
    \end{align*}
    It follows that $\dimAs\phi F \geq \alpha$, as required.
\end{proofref}
\begin{remark}\label{r:no-upper-box-dim-func}
    Since dimension functions $\phi(R)$ must decrease as $R$ decreases to $0$, there exists a dimension function such that $\dimAs\phi F = \dimuB F$ if and only if there is a number $\theta\in (0,1)$ so that $\dimAs\theta F = \dimuB F$.
\end{remark}
\subsection{Typicality of Moran sets}\label{ss:moran}
Throughout this section, we use the max norm on $\R^d$.

We give some motivation for the construction which we will use to prove the main result in this section.
Suppose we are given an infinite set $F\subset\R^d$ with diameter $1$ (with respect to the max norm).
For $n \in \N$, inductively define $\rho_n$ by
\begin{equation*}
    \rho_n \coloneqq \inf\{  r \in (0,1] : N_r(F) < 2^{(n+1)d}  \}.
\end{equation*}
Since $F$ is infinite, $\rho_n>0$.
Moreover, since a hypercube of side-length $r$ can be covered by $2^d$ hypercubes of side-length $r/2$, each $\rho_m \leq \rho_{m-1}/2$, and therefore the ratios $r_m\coloneqq\rho_{m}/\rho_{m-1}$ (using $\rho_0=1$) define a Moran set $M$.
We can then verify directly that $M$ satisfies
\begin{equation*}
    2^{-d} N_r(F) \leq N_r(M) \leq 2^d N_r(F).
\end{equation*}
In this situation, it follows immediately that $\dimuB F=\dimuB M$ and $\dimlB F=\dimlB M$.
In fact, since the Moran set has `average' branching everywhere, for all $0<r\leq R < 1$, it moreover holds that
\begin{equation*}
    N_r(M \cap B(x,R)) \leq 2^d \sup_{z\in F} N_r(F\cap B(z,R)).
\end{equation*}
In the actual proof, we will repeat this construction relative to a sequence of hypercubes $Q_n$ chosen to capture the worst-case scaling behaviour for a sequence of dimension functions.

In \cref{l:Moran-cnst}, we will formalize this construction.
We first require the following standard property of homogeneous Moran sets, which states that the covering numbers are approximately multiplicative and approximately constant at each scale.
For general sets, it only holds that maximal covering numbers are approximately sub-multiplicative.
\begin{lemma}\label{l:Moran-multi}
    Let $d\in\N$ be arbitrary.
    Then for all homogeneous Moran sets $M\subset\R^d$, all $0<\delta \leq r\leq R<1$ and all $x,y,z\in M$,
    \begin{equation}\label{e:equivs}
        N_\delta \bigl(M\cap B(x,R)\bigr) \approx_d N_\delta\bigl(M\cap B(y,r)\bigr)\cdot N_r\bigl(M\cap B(z, R)\bigr).
    \end{equation}
    Here the implicit constants depend only on the ambient dimension.
\end{lemma}
\begin{proof}
    Let $M$ be a homogeneous Moran set relative to the contraction ratios $(r_n)_{n=1}^\infty$.
    Let $0<r\leq R < 1$ be arbitrary, and for $t\in\{r, R\}$, let $m_t\in\N$ be such that $r_1\cdots r_{m_t} \leq t < r_1 \cdots r_{m_t-1}$, taking the empty product to be $1$.
    Then by the construction of the Moran set (see for instance the proof of \cref{p:moran-dim-formula}),
    \begin{equation*}
        N_r\bigl(M\cap B(x,R)\bigr) \approx_d 2^{d(m_r - m_R)}
    \end{equation*}
    for $x\in M$.
    Substituting this formula for the various expressions in \cref{e:equivs} completes the proof.
\end{proof}
Now we formalize the construction sketched at the beginning of this section as follows.
One can think of the conditions on $M$ below as \cref{im:lower} imitating the scaling of part of the set $F$ on the blocks $[\delta_n, R_n]$, without \cref{im:upper} being larger than $F$; and \cref{im:between} being as small as possible in the remaining ``gaps'' $[R_{n+1},\delta_n]$.
\begin{lemma}\label{l:Moran-cnst}
    Let $d\in\N$.
    There is a constant $C = C(d)\geq 1$ such that the following holds.
    Let $F\subset\R^d$ be arbitrary, let $1/2 \geq R_1 \geq \delta_1 \geq R_2 \geq \delta_2\geq \cdots > 0$ be a sequence of scales converging to $0$, and for each $n\in\N$ let $x_n\in F$.
    Then there exists a homogeneous Moran set $M\subset\R^d$ such that:
    \begin{enumerate}[nl,r]
        \item\label{im:lower} For all $n\in\N$ and $\delta_n \leq r \leq R_n$,
            \begin{equation*}
                \inf_{x\in M}N_r\bigl(M \cap B(x,R_n)\bigr) \geq C^{-1}\cdot N_r\bigl(F \cap B(x_n, R_n)\bigr).
            \end{equation*}
        \item\label{im:upper} For all $n\in\N$ and $\delta_n\leq r \leq R \leq R_n$,
            \begin{equation*}
                \sup_{x\in M} N_r\bigl(M \cap B(x,R)\bigr) \leq C\cdot \sup_{x\in F} N_r\bigl(F\cap B(x,R)\bigr).
            \end{equation*}
        \item\label{im:between} For all $n\in\N$,
            \begin{equation*}
                \sup_{x\in M} N_{R_{n+1}}(M\cap B(x, \delta_n)) \leq C.
            \end{equation*}
    \end{enumerate}
\end{lemma}
\begin{proof}
    First, assume that $\delta_n/R_n < 1/2$ for infinitely many $n$.
    If $\delta_n/R_n\geq 1/2$ for some $n\in\N$, then \cref{im:lower} holds with the trivial lower bound $N_{\delta_n}\bigl(M \cap B(x,R_n)\bigr) \geq 1$ by taking $C \geq 2^d$, regardless of the choice of $M$.
    Also, \cref{im:between} only becomes a stronger statement if there are fewer scales $\delta_n$ and $R_n$, and in particular implies \cref{im:upper} on the removed scales.
    Therefore we may assume that $\delta_n/R_n < 1/2$ for all $n\in\N$ by removing pairs of scales which do not satisfy this condition.

    We now proceed with the inductive construction of a sequence of contraction ratios $(r_n)_{n=1}^\infty$ which will define the homogeneous Moran set $M$.
    Begin with $r_1 = R_1$.

    Now suppose by induction that we have defined $r_1,\ldots,r_{j_n}$ and $\rho_1,\ldots,\rho_{j_n}$, related by $\rho_k = \prod_{j=1}^k r_j$, such that $R_n \geq \rho_{j_n} > R_n / 2$.
    For $m = j_n+1, j_n + 2,\ldots$, we inductively define $r_m'$ to be such that
    \begin{equation}\label{e:rm-choice}
        \rho_{j_n} \prod_{i=j_n+1}^m r_i' = \inf\{ r \in (0,R_n] : N_r(F \cap B(x_n, R_n)) \leq 2^d\cdot 2^{d(m - j_n)}  \},
    \end{equation}
    halting if $r_m'=0$ and setting all remaining terms equal to $0$.
    Let us verify by induction that $r_m' \leq 1/2$ for all $m\geq j_n + 1$.
    For the case $m = j_n+1$, since $N_{R_n/4}(B(x_n, R_n)) \leq 2^{2d}$ and $N_r(\cdot)$ is monotonic in $r$, we have $\rho_{j_n}r_{j_n+1}' \leq R_n/4$.
    Since $\rho_{j_n} \geq R_n/2$, therefore $r_{j_n+1}' \leq 1/2$.
    Now suppose $r_m'\in(0,1/2]$ for some $m \geq j_n+1$ and write $r'=\rho_{j_n} \prod_{i=j_n+1}^m r_i'$.
    By the definition of $r'$, for all $r > r'$, $N_r(F \cap B(x_n,R_n)) \leq 2^d \cdot 2^{d(m-j_n)}$.
    Therefore by covering each ball $N_r(B(x_n, R_n))$ by $2^d$ balls of radius $r/2$,
    \begin{equation*}
        N_{r/2}(F\cap B(x_n, R_n)) \leq 2^d \cdot 2^{d(m + 1-j_n)}.
    \end{equation*}
    Since $r>r'$ was arbitrary, it follows that $r_{m+1}' \leq 1/2$.

    We now choose the stopping index $j_{n+1} > j_n$ such that $\rho_{j_{n+1}-1} \approx \delta_{n}$.
    More precisely, let $k \geq 0$ be maximal such that (taking the empty product to be $1$)
    \begin{equation*}
        \rho_{j_n} \prod_{i=j_n+1}^{j_n+k} r_i' \geq \delta_n
    \end{equation*}
    and let $j_{n+1} = j_n + k + 1$.
    Then for $m\in \{j_n+1,\ldots,j_{n+1}-1\}$ (possibly there are no such $m$), let $r_m = r_m'$.
    Note that $\rho_{j_{n+1}-1} \geq \delta_n$ by the choice of $k$ above and since $\delta_n/R_n<1/2$, so we may choose $r_{j_{n+1}}$ such that $R_{n+1} \geq \rho_{j_{n+1}}\geq R_{n+1}/2$.
    Thus the induction may continue.
    Finally, let $M$ denote the homogeneous Moran set corresponding to the contraction ratios $(r_n)_{n=1}^\infty$.

    We now verify the desired properties of the construction.
    Let us first observe the following key consequences of \cref{e:rm-choice}.
    Let $n\in\N$ and let $\delta_n \leq r \leq R_n$ be arbitrary.
    Let $m\geq j_n$ be such that $2\rho_{m+1} < r \leq 2\rho_m$ (this choice is possible since $2\rho_{j_n} \geq R_n\geq r$).
    The choice of $m$ implies, since $r\approx \rho_m$ and $R_n \approx \rho_{j_n}$, that
    \begin{equation}\label{e:same-cov-init}
        2^{d(m-j_n)} \approx_d \inf_{x\in M}N_r\bigl(M\cap B(x,R_n)\bigr)\approx_d\sup_{x\in M}N_r\bigl(M\cap B(x,R_n)\bigr).
    \end{equation}
    We next verify that
    \begin{equation}\label{e:same-cov}
        N_r(F \cap B(x_n, R_n)) \approx_d 2^{d(m-j_n)}.
    \end{equation}
    First observe that $m \leq j_{n+1}$.
    If $m \leq j_{n+1}-1$, then \cref{e:rm-choice} combined with the fact that the covering number $N_r(\cdot)$ has discontinuities of size at most $2^d$ implies \cref{e:same-cov}.
    Otherwise, $m = j_{n+1}$.
    In this case, $\rho_{j_{n+1}-1} r_m' < \delta_n$ by definition of $j_{n+1}$, so since $r \geq \delta_n$, it follows from \cref{e:rm-choice} that $N_r(F\cap B(x_n, R_n)) \lesssim_d 2^{d(m-j_n)}$.
    But the covering number is monotonically increasing, and the lower bound holds for $j_{n+1}-1$, yielding \cref{e:same-cov}.

    Combining \cref{e:same-cov-init} and \cref{e:same-cov} shows \cref{im:lower}.

    Now to verify \cref{im:upper}, let $0<\delta_n \leq r\leq R \leq R_n <1$ be arbitrary.
    First, by covering balls of radius $R$ by with balls of radius $r$,
    \begin{equation*}
        N_r(F\cap B(x_n, R_n)) \leq N_R(F\cap B(x_n, R_n))\cdot \sup_{x\in F}N_r(F\cap B(x, R))
    \end{equation*}
    But the covering number of a homogeneous Moran set is approximately multiplicative as proven in \cref{l:Moran-multi}, so by \cref{e:same-cov-init} and \cref{e:same-cov} applied at scale $r$ and then $R$, for all $y\in M$,
    \begin{equation*}
        \sup_{x\in F}N_r(F\cap B(x, R)) \gtrsim_d \frac{N_r(M\cap B(x_n, R_n))}{N_R(M\cap B(x_n, R_n))} \approx_d N_r(M\cap B(y,R)).
    \end{equation*}
    This gives \cref{im:upper}, as required.

    Also, \cref{im:between} follows since $\rho_{j_{n+1}-1} \geq \delta_{n}$ and $\rho_{j_{n+1}}\leq R_{n+1}$ for all $n\in\N$, so there is only $1$ level between scales $\delta_n$ and $R_{n+1}$ in the construction of $M$.
    
    Finally, if $\delta_n/R_n\geq 1/2$ for all but finitely many $n$, the construction is much easier and is left to the interested reader.
\end{proof}
Using this technical lemma, we now obtain the following special (but key) case of \cref{it:morantypical}.
\begin{lemma}\label{l:moran-lemma}
    Fix $d \in \N$ and $F \subset \R^d$, and let $\{\phi_i\}_{i \in \N}$ be a countable family of dimension functions.
    Then there exists a homogeneous Moran set $M \subset \R^d$ such that $\dimAs{\phi_i} F \leq \dimAs{\phi_i}M$ for all $i \in \N$, and moreover $\dimAs{\psi} M \leq \dimAs\psi F$ for all dimension functions $\psi$.
    In particular, $\dimA M \leq \dimA F$.
\end{lemma}
\begin{proof}
    Fix an enumeration $(i_n)_{n=1}^\infty$ of $\N$ which contains each element of $\N$ infinitely often.
    Intending to use \cref{l:Moran-cnst}, we construct scales $R_n$ and $\delta_n$ inductively as follows.
    Set $\delta_0=1/2$, and inductively for $n\in\N$, get $x_{n}\in F$ and $0<R_{n}\leq \delta_{n-1}^{n}$ such that with $\delta_n = R_n^{1+\phi_{i_n}(R_n)}$,
    \begin{equation*}
        N_{\delta_n}(F\cap B(x_n, R_n)) \geq R_n^{-( \dimAs{\phi_{i_n}} F -  \frac{1}{n} )\phi_{i_n}(R_n)}.
    \end{equation*}
    Apply \cref{l:Moran-cnst} to the sequence of scales $R_1 \geq \delta_1 \geq R_2 \geq \delta_2 \geq \cdots$ to get a homogeneous Moran set $M$ satisfying the conclusions of \cref{l:Moran-cnst}.
    It is immediate from \cref{l:Moran-cnst}~\cref{im:lower} that $\dimAs{\phi_i} F \leq \dimAs{\phi_i}M$ for all $i\in\N$.

    For the other bound, let $\psi$ be an arbitrary dimension function and let $L\in\N$ be so that $L\geq \sup_{x\in (0,1)} \psi(x)$.
    Then let $0<r = R^{1+\psi(R)} \leq R < \delta_{L+1}$ be arbitrary.
    Since $R_n \leq \delta_{n-1}^n$ for all $n\in\N$, it follows from the choice of $L$ that there is an $m\in\N$ so that $R_{m+1} \leq r \leq R \leq \delta_{m-1}$.
    If $R\leq \delta_m$ or $R_m \leq r$, then it follows from \cref{l:Moran-cnst}~\cref{im:between} that $N_r(M\cap B(x,R)) \approx 1$.
    Otherwise, let $\delta_m \leq r' \leq R' \leq R_m$ be minimal (resp.\ maximal) such that $r \leq r' \leq R' \leq R$.
    Then applying \cref{l:Moran-cnst}~\cref{im:between} on the scales $[R_m,\delta_{m-1}]\cup [R_{m+1},\delta_m]$ followed by \cref{l:Moran-cnst}~\cref{im:upper} on the scales $[\delta_m,R_m]$, for any $\eps>0$,
    \begin{align*}
        N_r\bigl(M\cap B(x,R)\bigr)
        &\lesssim N_{r'}\bigl(M\cap B(x,R')\bigr)\\
        &\lesssim \sup_{x\in F}N_{r'}\bigl(F\cap B(x,R')\bigr)\\
        &\leq\sup_{x\in F}N_{r}\bigl(F\cap B(x,R)\bigr)\\
        &\lesssim_\eps R^{-\psi(R)(\dimAs\psi F +\eps)}.
    \end{align*}
    In either case, it follows that $\dimAs\psi M \leq\dimAs\psi F$.

    The ``in particular'' statement follows by taking a dimension function $\psi$ with $\dimAs\psi M = \dimA M$, which is guaranteed by \cref{t:Assouad-recover}.
\end{proof}
Finally, we can prove the main result.
For the convenience of the reader, we include a full statement of \cref{it:morantypical} here.
\begin{restatement}{it:morantypical}
    Fix $d \in \N$ and $F \subset \R^d$, and let $\mathcal{A}$ be a family of dimension functions.
    Suppose $\mathcal{A}=\bigcup_{i=1}^\infty\mathcal{A}_i$ where for each $i$ there exists $T_i \subset \R$ such that $\mathcal{A}_i = \{ \phi_{i,t} : t \in T_i\}$ and whenever $t,t' \in T_i$ satisfy $t \geq t'$ the following limit exists and lies in $[0,1]$:
    \begin{equation}\label{e:limit-cond}
        \lim_{R\to 0}\frac{\phi_{i,t}(R)}{\phi_{i,t'}(R)}\in[0,1].
    \end{equation}
    Then there exists a homogeneous Moran set $M \subset \R^d$ such that $\dimAs\psi M \leq\dimAs\psi F$ for all dimension functions $\psi$, and moreover
    \begin{equation*}
        \dimAs{\psi} F = \dimAs{\psi} M \qquad \text{and} \qquad \dimuAs{\psi} F = \dimuAs{\psi} M
    \end{equation*}
    for all $\phi \in \mathcal{A}$ and $\psi\in\mathcal{W}_{\phi}$.
\end{restatement}
\begin{proof}
    For $\phi\in\mathcal{A}$, write $s_\phi=\dimAs{\phi} F$.
    For all $i\in\N$, we choose a countable family of functions $\mathcal{C}_i\subset\mathcal{A}_i$ such that for all $\phi\in\mathcal{A}_i$ and $\eps>0$ there is $\psi\in\mathcal{C}_i$ such that
    \begin{equation}\label{e:sandwich}
        s_\phi - \varepsilon \leq s_{\psi}\qquad\text{and}\qquad \lim_{R\to 0}\frac{\phi(R)}{\psi(R)} < \infty.
    \end{equation}
    We choose such a family $\mathcal{C}_i$ as follows.
    Let $S_i = \{s_\phi:\phi\in\mathcal{A}_i\}\subset[0,d]$.
    First, let $\mathcal{A}_{i,0}$ denote a countable subset of $\mathcal{A}_i$ for which $\{s_\psi:\psi\in\mathcal{A}_{i,0}\}$ is dense in $S_i$.
    Next, for $m\in\N$, let $E_{i,m} = \{s_\phi\in S_i:(s_\phi - 1/m, s_\phi) \cap S_i = \varnothing\}$ and observe that $E_{i,m}$ is a finite set.
    Then for each $s\in E_{i,m}$, let $T_{i,s}$ denote the set of indices $t \in T_i$ so that $s_{\phi_{i,t}} = s$.
    If $T_{i,s}$ is bounded below and contains its infimum $t$, let $\mathcal{C}_{i,s}=\{\phi_{i,t}\}$; otherwise, let $\mathcal{C}_{i,s}=\{\phi_{i,t_n}\}$ where $t_n\in T_{i,s}$ converge to the infimum or diverges to minus infinity if $T_{i,m}$ is unbounded below.
    In any case the point is that if $t \in T_{i,s}$ then there is a $t'\leq t$ such that $\phi_{i,t'}\in\mathcal{C}_{i,s}$ and $s_{\phi_{i,t'}} = s$.
    Finally let
    \begin{equation*}
        \mathcal{C}_i = \mathcal{A}_{i,0} \cup \bigcup_{m=1}^\infty \bigcup_{s \in E_{i,m}}\mathcal{C}_{i,s}
    \end{equation*}
    which is a countable set.

    We now verify \cref{e:sandwich} for $\mathcal{C}_i$.
    Let $\phi\in\mathcal{A}_i \setminus \mathcal{C}_i$ be arbitrary and $\eps>0$.
    First, suppose there exists $\psi\in\mathcal{C}_i$ such that $s_\phi -\eps \leq s_\psi < s_\phi$.
    Since $s_\psi < s_\phi$, by \cref{it:fund-dim}~\cref{im:window-compare} applied with the set $F$, we must have $\lim_{R\to 0}\phi(R)/\psi(R) < \infty$.
    Otherwise, $s_\phi \in E_{i,m}$ for some $m$ and write $\phi = \phi_{i,t}$.
    By construction, get $\psi\in\mathcal{C}_{i, s_\phi}$ such that $\psi = \phi_{i, t'}$ where $t'\leq t$.
    Observe that $s_\psi = s_\phi$; and moreover, by \cref{e:limit-cond} since $t' \leq t$, $\lim_{R\to 0}\phi(R)/\psi(R)<\infty$ as required.

    Finally, let $\mathcal{C}$ denote the closure under multiplication by positive rationals of the union $\bigcup_{i=1}^\infty\mathcal{C}_i$.
    By \cref{l:moran-lemma}, get a homogeneous Moran set $M$ such that $\dimAs\phi M=\dimAs\phi F$ for all $\phi\in\mathcal{C}$, and moreover $\dimAs\psi M\leq \dimAs\psi F$ for all dimension functions $\psi$.

    Now suppose $i\in\N$ and $\phi\in\mathcal{A}_i$ is arbitrary.
    First, suppose
    \begin{equation*}
        \lim_{R\to 0}\frac{\psi(R)}{\phi(R)}\in(0,\infty)
    \end{equation*}
    for some $\psi\in\mathcal{C}_i$.
    Let $\eps>0$ and let $q\in\Q\cap(0,\infty)$ be such that
    \begin{equation*}
        \lim_{R\to 0}\frac{\phi(R)}{\psi_q(R)}\in(1-\eps,1+\eps).
    \end{equation*}
    Then by \cref{c:phi-equivs}, letting $M_d$ denote the doubling constant in $\R^d$,
    \begin{equation*}
        |\dimAs{\psi_q}M-\dimAs{\phi}M|\leq\eps(1+2\log_2 M_d+\eps),
    \end{equation*}
    and similarly for $F$ in place of $M$.
    But $\psi_q\in\mathcal{C}$ and $\eps>0$ was arbitrary, giving that $\dimAs\phi F = \dimAs\phi M$.

    Otherwise, for all $\psi\in\mathcal{C}_i$,
    \begin{equation*}
        \lim_{R\to 0}\frac{\psi(R)}{\phi(R)}\in\{0,\infty\}.
    \end{equation*}
    Let $\eps>0$ be arbitrary.
    Then by \cref{e:sandwich}, get $\psi\in\mathcal{C}_i$ so that $s_\phi-\eps\leq s_\psi$ and $\lim_{R\to 0}\phi(R)/\psi(R)=0$.
    In particular, by \cref{it:fund-dim}~\cref{im:window-compare} applied with the set $M$, it follows that
    \begin{equation*}
        s_\phi-\eps \leq s_\psi = \dimAs\psi F=\dimAs{\psi}M\leq\dimAs{\phi} M.
    \end{equation*}
    Since $\eps>0$ was arbitrary, it follows that $\dimAs\phi M \geq \dimAs\phi F$; and we recall that the other bound always holds, yielding the desired inequality.

    Finally, the result for the upper $\phi$-Assouad dimensions follows from \cref{it:upper-recover} and continuity \cref{it:fund-dim}~\cref{im:cont}, since $\mathcal{C}$ is closed under multiplication by any $q\in\Q\cap(0,\infty)$.
\end{proof}
\section{Stochastically self-similar sets}\label{s:stoch}
We now turn our attention to specific families of sets.
In this section, we focus on stochastically self-similar sets.
\subsection{Galton--Watson processes and random trees}\label{ss:gw-proc}
Let $X$ be a random variable taking values in the non-negative integers.
We say that $X$ is an \defn{offspring number} and refer to its distribution as an \defn{offspring distribution}.
Its associated probability generating function is
\begin{equation*}
  f(s) = \E(s^X) = \sum_{j=0}^{\infty}\Prob\{X=j\} \cdot s^j.
\end{equation*}
We say that $X$ is finitely supported if $N\coloneqq \max\left\{ n : \Prob\left\{ X=n \right\}>0
\right\}$ is finite.
The associated probability generating function is then a polynomial of degree $N$ with non-negative coefficients $\theta_j = \Prob\left\{ X=j \right\}$ such that $\sum_{i=0}^N \theta_j =1$.

We recall some basic properties of $f$, a proof for which can be found in \cite{zbl:0259.60002}.
\begin{proposition}\label{p:pgf}
    Let $f$ be the probability generating function of a non-negative integer valued random variable $X$ such that $\E(X)>0$.
    Then:
    \begin{enumerate}[nl,r]
        \item $f$ is smooth, convex, and strictly increasing on $[0,\infty)$.
            It is strictly convex if and only if there exists $n\geq 2$ such that $\Prob\left\{ X=n \right\}>0$ (the non-trivial case).
        \item $f(0) = \theta_0$ and $f(1) = 1$.
        \item The expectation of $X$ is $m\coloneqq\E(X) = f'(1)$.
        \item If $m=\E(X)>1$ (the supercritical case), then $f$ is non-trivial and there exists a unique $q\in [0,1)$ such that $f(q) = q$.
    \end{enumerate}
\end{proposition}

The Galton--Watson process $Z_k$ with offspring variable $X$ is defined by the recursion
\begin{equation*}
  Z_0=1\quad
  \text{and}
  \quad
  Z_{k+1} = \sum_{i=1}^{Z_k} X_{k,i}\,,
\end{equation*}
where $X_{k,i}$ are independent random variables that equal $X$ in distribution.
We can relate the probability generating function of $X$ to the behaviour of the Galton--Watson process, see \cite{zbl:0259.60002} for details.

\begin{proposition}
    Let $Z_k$ be a supercritical Galton--Watson process with offspring random variable $X$.
    Denote the probability generating function of $X$ by $f$.
    Then,
    \begin{enumerate}[nl,r]
        \item The probability generating function of $Z_k$ is $\E(s^{Z_k}) = f_k(s)$, where $f_k$ is the $k$-fold composition of $f$.
        \item The mean of $Z_k$ is given by $\E(Z_k) = \tfrac{d}{ds}f_k(s)\rvert_{s=1} = f_k'(1) = m^k$.
        \item The process $Z_k$ dies out with probability $q$, i.e.\ $\Prob(Z_k =0 \text{ for some }k)=q$ where $q\in[0,1)$ is the unique number such that $f(q)=q$.
    \end{enumerate}
\end{proposition}

From now we will assume that $X$ is a non-trivial finitely supported offspring random
variable. In this section we determine the Assouad spectrum for arbitrary dimension functions of the
Gromov boundary of its associated tree. 

We first define the Galton--Watson tree.
Let $\Lambda = \left\{ 1,\dots,N \right\}$ be a finite alphabet of size $N$, i.e.\ the degree of the probability generating function $f$.
Let $\Lambda^k$ denote the set of words of length $k$ over the alphabet $\Lambda$, and we let $\Lambda^*=\bigcup_{j=0}^\infty \Lambda^j$ denote the set of all finite words where $\Lambda^0 = \left\{ \varnothing \right\}$ contains only the empty word.
We let $\Lambda^{\N}$ denote the set of infinite words over $\Lambda$.

For $v\in\Lambda^*$ let $X_v$ be a random variable with the same distribution as $X$, independent of all distinct words $w\in\Lambda^*$.
The random Galton--Watson tree $\cT=\cT(\omega)$ is defined inductively by the rules
\begin{equation*}
  \cT_0 = \left\{ \varnothing \right\},\qquad
  L_0 = \cT_0,\qquad
  L_{n+1} = \bigcup_{v\in L_n}\left\{ vj\in\Lambda^{n+1} : 1\leq j \leq X_v \right\},
\end{equation*}
and
\begin{equation*}
  \cT_{n+1} = \cT_{n} \cup L_{n+1},
  \qquad
  \cT = \bigcup_{n=1}^{\infty}\cT_n.
\end{equation*}
Its boundary $\partial\cT$ is
\begin{equation*}
  \partial \cT = \left\{ v\in\Lambda^{\N} : v\rvert_n \in \cT \text{ for all }n\in\N \right\}
\end{equation*}
which consists of all ``eventually surviving'' branches of the finite trees.
Note that almost surely $\partial\cT$ is either the empty (extinct) tree or an infinite subtree of the full $N$-ary tree $\Lambda^{\N}$.
We define a metric on $\Lambda^{\N}$, and thus $\partial\cT$, by $d(v,w) = e^{-|v\wedge w|}$, where $v\wedge w$ is the longest common ancestor of $v$ and $w$.
Note that $v\wedge w=\varnothing$ and $v\wedge w = v \Leftrightarrow v=w$ are possible outcomes.
The metric space $(\partial\cT(\omega), d)$ is known as the Gromov boundary of the tree
$\cT(\omega)$.
Since the offspring distribution is assumed to be finitely supported, every ball of radius
$r=e^{-k}$ can be covered by $N$ balls of radius $r/e=e^{-(k+1)}$. Hence, $\partial\cT$ is doubling.

\subsection{Large deviations of Galton--Watson processes}\label{s:large-deviations}
Define $\gamma$ such that $m^\gamma = N$.
Note that $\gamma\geq 1$ with equality only occurring when $m=N$, i.e.\ when $X = N$ is constant almost surely.
In this section, we will prove \cref{p:GW-large-deviations}, which we reproduce here for the convenience of the reader.
\begin{restatement}{p:GW-large-deviations}
Let $Z_k$ be a Galton--Watson process with offspring random variable $X$ which is not almost surely constant.
    Assume that its probability generating function $f$ is a polynomial of degree $2\leq N <\infty$ and $m\coloneqq \E(X)>1$.
    Define $\gamma$ such that $m^\gamma = N$.
    Then for all $1<t<\gamma$, all $\eps>0$ sufficiently small, and all $k \in \N$,
    \begin{equation*}
        \exp\left( -m^{(t-1+\eps)\frac{\gamma}{\gamma-1}k} \right)\lesssim
        \Prob\left( Z_k  \geq m^{tk} \right)
        \lesssim \exp\left( -m^{(t-1-\eps)\frac{\gamma}{\gamma-1}k} \right),
    \end{equation*}
    with the implicit constants depending only on $t$ and $\eps$.
\end{restatement}
Before proceeding with a proof of this result, we remark that similar results have been known for a long time.
However, they often concern estimates on the limiting variable $W \coloneqq \lim_k Z_k / m^k$, which exists almost surely.
In particular, Harris \cite{zbl:0041.45603} showed that
\begin{equation*}
    \log\E\left( e^{s W} \right) = s^\gamma H(s)+ O(1),
\end{equation*}
where $H$ is a continuous, positive, and multiplicatively periodic function (cf.\ \cref{thm:polybound} below).
The probabilistic analogue was derived by Biggins \& Bingham \cite{zbl:0796.60090}, who showed that
\begin{equation*}
    \log\Prob\left( W > x \right) = x^{\gamma/(\gamma-1)}\widetilde{H}(x)+ o(x^{\gamma/(\gamma-1)}).
\end{equation*}
Further extensions have been proved in \cite{zbl:1032.60048}.
While the behaviour of $W$ can be related back to that of $W_k\coloneqq Z_k/m^k$, see e.g.~\cite{zbl:0806.60068}, we needed explicit bounds on the rate of convergence and give a self-contained account here.
Further results are known on the tail behaviour of $W$, even if $X$ is not finitely supported, or even heavy-tailed, see for example \cite{zbl:1175.60075} and \cite{zbl:1290.60087}.
We point out \cite[Theorems 3 \& 4]{zbl:1290.60087} in particular, which characterizes the behaviour of $\Prob\left( W_k>x \right)$ for more general distributions.

It is plausible that \cref{p:GW-large-deviations} can be similarly sharpened through the application of a suitable renewal theorem.
However, the bounds in \cref{p:GW-large-deviations} are sufficient to establish our Borel--Cantelli lemma for trees (\cref{l:tree-bc}) so we have not attempted this here.

Our large deviations result depends fundamentally on the following asymptotic result for polynomial functions.
Recall that $f_k$ denotes the $k$-fold composition of the probability generating function $f$.
\begin{lemma}\label{thm:polybound}
    Let $f$ be a probability generating function of degree $2 \leq N <\infty$.
    Let $\eps>0$ and assume $m = f'(1)>1$ and $m<N$.
    Then for all $s\in (0,1)$,
    \begin{equation*}
        \lim_{k\to\infty} \frac{1}{N^{(s+\eps)k}}\log f_{k}\left( \exp\left( m^{-(1-s)k} \right) \right) = 0.
    \end{equation*}
\end{lemma}
\begin{proof}
    Let $\eps>0$ and fix $\delta>0$ such that $(1+\delta)^2\leq m^{\eps/2}$.
    Let $x_\delta>1$ be such that $f'(x_\delta) = (1+\delta)m \eqqcolon m_{\delta}$.
    Then, for $1\leq y \leq x_\delta$,
    \begin{equation*}
        f(y) \leq m_{\delta} (y-1)+1.
    \end{equation*}
    Now
    \begin{equation*}
        y_k \coloneqq \exp(m^{-(1-s)k}) = 1+m^{-(1-s)k}+O(m^{-2(1-s)k})< 1+(1+\delta)m^{-(1-s)k}
    \end{equation*}
    for all sufficiently large $k$.
    Let $n_k = \min\left\{k,\min\left\{ j : f_j(y_k) > x_\delta\right\}\right\}$.
    We can obtain a lower bound on $n_k$ for all sufficiently large $k$ by noting that
    \begin{equation*}
        f_{n_k}(y_k) \leq m_{\delta}^{n_k} (y_k-1)+1 \leq 1+(1+\delta)^{n_k+1}m^{n_k-(1-s)k} \leq
	1+m^{n_k-(1-s-\eps/2)k},
    \end{equation*}
    which gives $n_k\geq (1-s-\eps/2)k+c_\delta$, where $c_\delta \coloneqq \log_m(x_\delta-1)$.
    Moreover, since $f$ is a polynomial of degree $N$ with positive coefficients that are bounded above by $1$, we also have the trivial bound $f(x) \leq 1+(x-1)^N$ for all $x\geq 1$.
    Combining this with the previous estimate,
    \begin{align*}
        \log f_k(y_k) &= \log f_{k-n_k+1}\circ f_{n_k-1}(y_k) \\
        &\leq \log\left(1+ x_\delta^{N^{k-n_k+1}}\right) \\
                      &\leq N^{k-n_k+1}\log x_\delta+\log 2\\
                      &\leq N^{(s+\eps/2)k+1-c_{\delta}}\log x_{\delta}+\log 2.
    \end{align*}
    The conclusion now follows after dividing by $N^{(s+\eps)k}$ and taking limits.
\end{proof}
The following probability estimate does not require independence; its
short and elementary proof can be found, for example, in \cite[Lemma~2.1]{zbl:1512.28008} and we include it 
for convenience.
\begin{lemma}\label{l:boundLow}
    Let $E_1,\dots,E_n$ be a sequence of events with $\Prob(E_i)\geq p$ for all $1\leq i \leq n$.
    Let $0<\lambda < p$ and let
    $F_\lambda$ denote the event that at least $\lambda n$ of the events $E_1,\dots,E_n$ occur.
    Then
    \begin{equation*}
        \Prob ( F_{\lambda} )\geq \frac{p-\lambda}{1-\lambda}.
    \end{equation*}
\end{lemma}
\begin{proof}
    We compute
    \begin{align*}
        p n &\leq \E(\#\left\{ i : E_i \text{ occurs} \right\})\\
            &=\E(\#\left\{ i : E_i \text{ occurs} \right\} \mid F_\lambda)\Prob(F_\lambda) +
            \E(\#\left\{ i : E_i \text{ occurs} \right\} \mid F_{\lambda}^c)(1-\Prob(F_\lambda))
            \\
            &\leq n\Prob(F_\lambda) + \lambda n(1-\Prob(F_\lambda)),
    \end{align*}
    so that $p \leq \Prob(F_\lambda) +\lambda(1-\Prob(F_\lambda))$.
\end{proof}
We are now ready to complete the proof of \cref{p:GW-large-deviations}
\begin{proofref}{p:GW-large-deviations}
    We first prove the upper bound.
    Fix
    \begin{equation*}
        s = \frac{t-1}{\gamma -1} - \eps \frac{\gamma}{\gamma -1}
    \end{equation*}
    and note that
    \begin{equation}\label{eq:domination}
        s+t-1>(s+\eps/2)\gamma.
    \end{equation}
    Note also that $s>0$ and $s+\eps/2<1$ for $\eps>0$ small enough.
    The proof of the proposition now follows from Markov's inequality and a Chebyshev argument, as well as applying \cref{thm:polybound} with $\eps/2$.
    Indeed, there exists $C'>0$ such that for all sufficiently large $k$,
    \begin{align*}
        \Prob\left( Z_k  \geq m^{tk} \right)
    &= \Prob\left( \exp\left(m^{-(1-s)k} Z_{k}\right) \geq e^{m^{(s+t-1)k}} \right)\\
    &\leq \E\left( \exp\left(m^{-(1-s)k} Z_{k}\right)\right) \cdot e^{-m^{(s+t-1)k}}\\
    &= f_k\left( \exp\left( m^{-(1-s)k} \right) \right) \cdot e^{-m^{(s+t-1)k}}\\
    &\leq \exp\left( C' N^{(s+\eps/2)k} -m^{(s+t-1)k} \right)\\
    &=\exp\left( C' m^{(s+\eps/2)\gamma k}-m^{(s+t-1)k} \right).
    \end{align*}
    Note that by \cref{eq:domination} the negative term will eventually dominate, so we conclude that
    \begin{equation*}
        \Prob\left( Z_k  \geq m^{tk} \right) \leq C \exp\left( -m^{(s+t-1)k} \right) = C \exp\left( -m^{(t-1-\eps)\frac{\gamma}{\gamma-1}k} \right)
    \end{equation*}
    for some $C>0$ and all $k \in \N$, establishing the upper bound.

    We now prove the lower bound.
    First recall that the normalized Galton--Watson process $W_k = Z_k/m^k$ converges almost surely to a random variable $W$ with expectations satisfying $\E(W_k)=\E(W)=1$ for all $k\in\N$.
    Further, $\Prob(W>0) = 1-q$, where $q$ is the least root of $f(q)=q$.
    Hence there exists $\beta>0$ such that $\Prob(W_k>\tfrac12)>\beta$ for all $k$.
    Write $\theta_N = \Prob(X=N)$.
    The probability that $Z_1 = N, Z_2 = N^2, \dots, Z_n = N^n$ is
    \begin{equation*}
        \theta_N^{N+N^2+\dots+N^n} \geq \theta_N^{N^{n+1}}.
    \end{equation*}
    Let $\lambda = \tfrac12\beta$.
    Note that the probability that a Galton--Watson process satisfies $Z_{k-n}^{(i)} \geq \tfrac12 m^{k-n}$ can be bounded below by $\beta$.
    Thus, by \cref{l:boundLow}, the probability that out of $N^n$ independent realizations of $Z_{k-n}^{(i)}$ $(1\leq i
    \leq N^n)$ at least $\lambda\cdot N^n$ are larger than $\tfrac12m^{k-n}$ is bounded below by
    $\alpha = \tfrac{\beta}{2-\beta}>0$.
    It follows that
    \begin{equation*}
        \Prob\left(Z_{k}\geq
        \lambda N^n \tfrac12 m^{k-n}\right)
        =
        \Prob\left(Z_{k}
        \geq \tfrac12\lambda m^{\gamma n}m^{k-n}\right)
        \geq \alpha\;\theta_N^{N^{n+1}}
        =\alpha e^{-m^{\gamma (n+1)}\log\theta_N^{-1}}.
    \end{equation*}
    Letting $n$ be the least integer such that $m^{tk}\leq \tfrac12 \lambda m^{(\gamma-1) n+k}$
    gives
    \begin{equation*}
      n\leq \frac{t-1}{\gamma-1}k-\frac{\log_m(\lambda/2)}{\gamma-1}+1.
    \end{equation*}
    Therefore
    \begin{align*}
      \Prob\left( Z_k \geq m^{tk} \right)
      &\geq
      \Prob \left( Z_k\geq \tfrac12\lambda m^{(\gamma-1)n+k} \right)
      \\
      &\geq
      \alpha\exp\left(
	-m^{\frac{\gamma}{\gamma-1}(t-1)k}m^{-\frac{\gamma}{\gamma-1}\log_m(\lambda/2)+2\gamma}
      \log\theta^{-1}_N \right)
      \\
      &\geq
      \alpha \exp\left( -m^{\frac{\gamma}{\gamma-1}(t-1+\eps)k} \right)
    \end{align*}
    for all $k$ large enough such that
    \begin{equation*}
      m^{\frac{\gamma}{\gamma-1}\eps k}
      \geq
      m^{\frac{\gamma}{\gamma-1}\log_m(2/\lambda)+2\gamma} \log\theta_N^{-1}.
    \end{equation*}
    Our claim immediately follows.
\end{proofref}
To conclude this section, we also note that the lower bound can be improved to guarantee that a large number of the children also survive.
Recall that $q$ is the least root of $f(q)=q$.
\begin{corollary}\label{c:GW-deviations-surviving}
    Let $Z_k$ be a Galton--Watson process with offspring random variable $X$.
    Assume that its probability generating function $f$ is a polynomial of degree $2\leq N <\infty$ and $m=\E(X)>1$.
    Then for all $1<t<\gamma$, all $\eps>0$ sufficiently small, and all $k \in \N$,
    \begin{equation*}
        \Prob\left( Z_k  \geq m^{tk}\text{ and at least $\tfrac{1-q}{2}m^{tk}$ offspring processes survive}\right)
        \gtrsim\exp\left( -m^{(t-1+\eps)\frac{\gamma}{\gamma-1}k} \right)
    \end{equation*}
    with implicit constants depending only on $t$ and $\eps$.
\end{corollary}
\begin{proof}
    Note that all of the $m^{tk}$ offspring processes are independent of the event $\left\{ Z_k \geq m^{t k}\right\}$.
    The probability of survival is $1-q \in (0,1]$.
    Hence, by \cref{l:boundLow},
    \begin{align*}
        &\Prob\left( Z_k  \geq m^{tk} \text{ and more than $\tfrac{1-q}{2}m^{tk}$ offspring processes survive} \right)
        \\
        &=
        \Prob\left( Z_k \geq m^{tk} \right) \Prob\left( \text{at least $\tfrac{1-q}{2} m^{tk}$ independent processes survive} \right)
        \\
        &\geq
        \Prob\left( Z_k \geq m^{tk} \right)
        \frac{1-q}{1+q}.
    \end{align*}
    From this the claim immediately follows from \cref{p:GW-large-deviations}.
\end{proof}

\subsection{Borel--Cantelli for trees}
In this section we prove the useful Borel--Cantelli lemma for Galton--Watson trees mentioned in the introduction, namely \cref{l:tree-bc}.
\begin{restatement}{l:tree-bc}
  Let $E_k$ be any measurable event for a Galton--Watson tree and write $P_k = \Prob(E_k)$.
  Let $\widetilde{E}$ be the event that there are infinitely many $k\in\N$ such that a Galton--Watson
  tree contains a subtree $\cT(v)\in E_k$ at level $k$.
  \begin{enumerate}[nl,r]
    \item\label{im:BC0} $\Prob(\widetilde{E}) =0$ \;if\; $\sum_{n\in\N}P_n m^n <\infty$,
    \item\label{im:BC1} $\Prob(\widetilde{E}) =1$, conditioned on non-extinction, if there exists a summable sequence $K_n$ of non-negative numbers such that $\sum_{n\in\N}K_n P_n m^n=\infty$.
  \end{enumerate}
\end{restatement}

We remark that this lemma does not require an explicit independence condition that is stipulated for the second part of the standard Borel--Cantelli lemma.
Such an independence condition is replaced by the existence of a decreasing sequence $K_n$, which will allow us to use independent subtrees instead.
\begin{proof}
    We prove \cref{im:BC0} and \cref{im:BC1} separately.
    \begin{proofpart}
        Proof of \cref{im:BC0}.
    \end{proofpart}
    Assume first that $\sum_n P_n m^n<\infty$.
    Let $v\in L_k$ be a node in $\cT$ at level $k$.
    We write
    \(
    E_k(v) = \{\cT(v) \in E_k\}
    \)
    for the event that a given subtree $\cT(v)$ is in $E_k$.
    Similarly, for any subset $A\subset L_k$, we write
    \begin{equation*}
        E_k(A) = \{\exists v\in A \text{ such that }\cT(v) \in E_k\}.
    \end{equation*}
    Further, we let
    \begin{equation*}
        \widetilde{E} = \left\{ E_k(L_k) \text{ holds for infinitely many }k\in\N \right\}.
    \end{equation*}

    Note that $\#L_k/m^k$ is bounded almost surely.
    Therefore
    \begin{equation*}
        \widetilde{E} =
        \bigcup_{Q=1}^{\infty}\left(\{\limsup_{k\to\infty}E_k(L_k)\}\cap\{\#L_k \leq Q m^k\}\right)
    \end{equation*}
    and by continuity from below,
    \begin{equation*}
        \Prob( \widetilde{E} )
        =\lim_{Q\to\infty} \Prob\left( \left\{ \limsup_{k\to\infty}E_{k}(L_k)\right\} \bigm| \#L_k \leq Q m^k \right)
        \cdot\Prob\left( \#L_k\leq Q m^{k} \right).
    \end{equation*}
    Now $\Prob\left( \#L_k\leq Q m^{k} \right)$ increases to $1$ in $Q$, and by
    continuity from above,
    \begin{align*}
        \Prob\left( \left\{ \limsup_{k\to\infty}E_{k}(L_k)\right\} \bigm| \#L_k \leq Q m^k \right)
    &=\Prob\left( \bigcap_{n=1}^{\infty}\bigcup_{k=n}^{\infty} E_{k}(L_{k}) \bigm| \#L_k \leq Q m^k \right)
    \\
    &=
    \lim_{n\to\infty}\Prob\left( \bigcup_{k=n}^{\infty}
    E_{k}(L_{k})\mid \#L_k \leq Q m^k \right)
    \\
    &\leq \lim_{n\to\infty}\sum_{k=n}^\infty
    \Prob(E_k(L_k)\mid\#L_k \leq Q m^k).
    \end{align*}
    Using a standard Taylor bound,
    \begin{equation*}
        \Prob(E_k(L_k)\mid\#L_k \leq Q m^k) \leq 1-(1-P_k)^{Q m^k} \leq Q m^k P_k,
    \end{equation*}
    and so, writing $S_n = \sum_{k=n}^\infty P_km^k$,
    we obtain
    \begin{equation*}
        \Prob(\widetilde{E})
        \leq \lim_{Q\to\infty} \lim_{n\to\infty}\sum_{k=n}^\infty Q P_k m^k
	=\lim_{Q\to\infty}\left(Q\cdot\lim_{n\to\infty}S_n\right)=0
    \end{equation*}
    since $S_1<\infty$ by assumption.

    \begin{proofpart}
        Proof of \cref{im:BC1}.
    \end{proofpart}
    Assume now that $\sum_n K_n P_n m^n = \infty$.
    A standard calculation shows that $\prod_{n=2}^\infty (1-n^{-2}) = 1/2$.
    Clearly, $\sum_n \frac{1}{2} K_n P_n m^n = \infty$.
    Since $\sum_n K_n < \infty$, we have $K_n \to 0$ and we may assume without loss of generality that
    $K_n m^n < \tfrac1N m^n$ and that $K_n m^n$ is a non-negative integer by some bounded rescaling of $K_n$.

    Recalling that $W \coloneqq \lim_k Z_k / m^k$, note that
    \begin{align}
        \Prob\left( \widetilde{E} \right)
    &=\Prob\left( \widetilde{E} \mid W=0 \right) \Prob\left( W=0 \right)
    + \Prob\left( \widetilde{E} \mid W>0 \right)\Prob\left( W>0 \right)
    \nonumber
    \\*
    & =
    \Prob\left( \widetilde{E} \mid W>0 \right)(1-q).\label{eq:probConditioned}
    \end{align}
    The conditioning on $W>0$ implies that there almost surely exists a first split, i.e.\ a $k_0\in
    \N$ such that $\#L_{k_0-1}=1$ and $\#L_{k_0} >1$.
    Write $\tau_j$ for the probability that $\#L_{k_0}=j$, conditioned on non-extinction.
    By the independence of Galton--Watson processes we must have
    \begin{equation*}
        p_E \coloneqq \Prob\left( \widetilde{E}^c \mid W>0 \right)
        =\sum_{j=2}^N \tau_j (p_E)^j, \quad\text{where}\quad
        \sum_{j=2}^N \tau_j = 1.
    \end{equation*}
    A standard convexity argument now implies that the only real solutions in $[0,1]$ are $0$ and $1$.
    Hence $\Prob( \widetilde{E} \mid W>0 )\in\left\{ 0,1 \right\}$ and to show that $\widetilde{E}$
    occurs almost surely conditioned on the Galton--Watson process not going extinct, we need only to
    prove that $\Prob(\widetilde{E})>0$ by \cref{eq:probConditioned}.

    Let $k\in\N$ and $A\subseteq L_k$ and write $P_k = \Prob(E_k)$.
    By independence,
    \begin{align*}
        \Prob(E_k(A)) &= 1- \prod_{v\in A} \Prob(E_k(v)^c) = 1-\left( 1-P_k \right)^{\#A}.
    \end{align*}
    Write $a_k = 1-(k+1)^{-2}$ and note that $\prod_{k=1}^{\infty}a_k = 1/2$.
    Let $k_0\in\mathbb{N}$ be large enough such that $(N/m)^{k_0} > 2\sum_{j}K_j$.
    Recalling that $\theta_N = \Prob(X=N)$, with probability $(\theta_N)^{1+N^2+\cdots+N^{k_0-1}}>0$ we have $\#L_{k_0} = N^{k_0}$, and we shall call this event $F_{k_0}$.
    Let $A_{k_0} \subset L_{k_0}$ be such that $\#A_{k_0} = K_{k_0} m^{k_0}(< N^{k_0})$.
    For definiteness, and to avoid dependency on the event $F_{k_0}$, we define $A_{k_0}$ to be independently randomly chosen elements of $L_{k_0}$ stopping once $\#A_{k_0}$ reaches $\min\{\#L_{k_0}, K_{k_0} m^{k_0}\}$.
    The probability of $E_{k_0}(A_{k_0})$, conditioned on $\#L_{k_0}\geq K_{k_0}m^{k_0}$, is given by $\Prob(E_{k_0}(A_{k_0})) = 1-(1-P_{k_0})^{K_{k_0} m^{k_0}}$.
    For convenience we define
    \begin{equation*}
        M_{k_0} = N^{k_0} \quad\text{and}\quad
        M_{k+1} = (M_k - K_k m^k)\cdot a_{k} m
    \end{equation*}
    for $k \geq k_0$.
    It is straightforward to see that
    \begin{equation*}
        M_k = m^k \left( \frac{N^{k_0}}{m^{k_0}} \prod_{j=k_0}^{k-1}a_j - \sum_{j=k_0}^{k-1} K_j \cdot\prod_{i=j}^{k-1}a_i \right)
        \geq m^{k}\left( \frac{N^{k_0}}{m^{k_0}}\prod_{j=1}^{\infty}a_j - \sum_{j=1}^{\infty}K_j \right)
        \geq m^k.
    \end{equation*}
    Consider now the descendants of $L_{k_0}\setminus A_{k_0}$.
    Conditioned on $F_{k_0}$,
    \begin{equation*}
        \#(L_{k_0}\setminus A_{k_0}) \geq N^{k_0} - K_{k_0}m^{k_0} \geq m^{k_0}.
    \end{equation*}
    Their offspring number is independent of $E_{k_0}(A_{k_0})$ and we write
    \begin{equation*}
        F_{k_0+1} = \left\{\# \left\{v\in L_{k_0+1} :  v|_{k_0} \in L_{k_0}\setminus A_{k_0}\right\}
        \geq M_{k_0+1}  \right\}.
    \end{equation*}
    for the event that its offspring number is at least $M_{k_0+1}$.

    Inductively, let $A_{k}\subseteq L_{k}$ be an arbitrary, randomly chosen, subset of cardinality equal to $\min\{\#L_k,K_km^k\}$.
    We can define $F_k$, for $k > k_0$, to be the event
    \begin{equation*}
        F_{k} = \left\{\# \left\{v\in L_{k} :  v|_{k-1} \in L_{k-1}\setminus A_{k-1}\right\}
        \geq M_{k}  \right\}.
    \end{equation*}
    We can estimate $\Prob\left( F_{k+1} \mid \langle F_{k_0},F_{k_0+1},\ldots,F_{k}\rangle \right)$ using Hoeffding's inequality:
    \begin{align*}
        \Prob&\left( F_{k+1} \mid \langle F_{k_0},F_{k_0+1},\ldots,F_{k}\rangle \right)\\*
             &\geq \Prob\left( \sum_{i=1}^{M_{k}-K_{k}m^{k}} X_{k,i} \geq M_{k+1} \right)
             \\
             &=\Prob\left( \sum_{i=1}^{M_{k}-K_{k}m^{k}} X_{k,i} \geq \left(M_{k}-K_km^k\right)a_k m \right)
             \\
             &=\Prob\left( -\sum_{i=1}^{M_{k}-K_{k}m^{k}} (X_{k,i}-m) \leq (1-a_k)\left(M_{k}-K_km^k\right) m \right)
             \\
             &=1-\Prob\left( -\sum_{i=1}^{M_{k}-K_{k}m^{k}} (X_{k,i}-m) > \frac{m}{(k+1)^2}\left(M_{k}-K_km^k\right)  \right)
             \\
             &\geq 1-\exp\left( -\,\frac{m^2}{(k+1)^4}\cdot\frac{\left(M_{k}-K_km^k\right)^2}{\left(M_{k}-K_km^k\right)\cdot N^2} \right)
             \\
             &= 1-\exp\left( -\,\frac{m^2}{N^2}\cdot\frac{M_{k}-K_km^k}{(k+1)^4} \right)
             \\
             &\geq 1-\exp\left( -\,\frac{m^2(N-1)}{N^3} \cdot \frac{m^k}{(k+1)^4} \right).
    \end{align*}
    In particular, this gives
    \begin{align*}
        \Prob&\left( F_{k_0}\cap \cdots\cap F_{k+1} \right)\\*
    &=\Prob\left( F_{k+1} \bigm| \langle F_{k_0},\ldots,F_{k} \rangle\right)\cdot
    \Prob\left( F_k \bigm| \langle F_{k_0},\ldots,F_{k-1} \rangle\right)
    \cdots
    \Prob\left( F_{k_0} \right)
    \\
    &\geq
    \left( \theta_N \right)^{1+N^2+\cdots+N^{k_0-1}}\cdot
    \prod_{j=k_0}^k
    \left(1-\exp\left( -\,\frac{m^2(N-1)}{N^3} \cdot \frac{m^j}{(j+1)^4} \right)\right)
    \\
    &\geq
    \left( \theta_N \right)^{1+N^2+\cdots+N^{k_0-1}}\cdot
    \prod_{j=k_0}^\infty
    \left(1-\exp\left( -\frac{m^j}{N^3(j+1)^4} \right)\right)\\
    &\eqqcolon \alpha >0,
    \end{align*}
    noting that the lower bound is independent of $k$.
    Writing $F=\bigcap_{k=k_0}^{\infty} F_k$, we get $\Prob(F) \geq \alpha$.

    Now
    \begin{equation*}
        \Prob(\widetilde{E})\geq \Prob(\widetilde{E}\cap F) = \Prob(\widetilde{E}\mid F)\Prob(F) \geq \alpha\cdot \Prob(\widetilde{E}\mid F).
    \end{equation*}
    Since $\widetilde{E} \supseteq \limsup_{k\to\infty}E_k(A_k)$, we also have
    \begin{equation*}
        \Prob(\widetilde{E}\mid F)\geq 1-\Prob\left( \bigcup_{k=k_0}^\infty \bigcap_{n=k}^\infty E_n(A_n)^c \bigm| F \right).
    \end{equation*}
    The events $E_n(A_n)$ and $E_{n'}(A_{n'})$ are independent for $n\neq n'$ since the vertices in
    $A_n$ are, by definition, not descendants of $A_{n'}$ and vice versa. Hence,
    \begin{align*}
        \Prob\left( \bigcup_{k=k_0}^\infty \bigcap_{n=k}^\infty E_n(A_n)^c \bigm| F \right)
    &\leq
    \sum_{k=k_0}^{\infty}
    \Prob\left( \bigcap_{n=k}^{\infty} E_{n}(A_{n})^c \bigm| F \right)
    \\
    &=
    \sum_{k=k_0}^{\infty}
    \prod_{n=k}^{\infty}\Prob\left( E_{n}(A_{n})^c \bigm| F \right)
    \\
    &=
    \sum_{k=k_0}^{\infty}
    \left( \lim_{K \to \infty} \exp\left(
        \sum_{n=k}^{K}\log\left( 1-P_n \right)^{K_n m^n} \right)
    \right)
    \\
    &\leq
    \sum_{k=k_0}^{\infty}
    \left(\lim_{K \to \infty}\exp\left(-
        \sum_{n=k}^{K}P_n K_n m^n \right)
    \right)
    \\
    &=0
    \end{align*}
    by our divergence assumption.
    We conclude that $\Prob(\widetilde{E}\mid F) = 1$ and so $\Prob(\widetilde{E}) \geq \alpha>0$.
    This completes the proof.
\end{proof}

\subsection{The \texorpdfstring{$\phi$}{𝜙}-Assouad dimensions of branching processes}
We now prove the exact formula for the $\phi$-Assouad dimensions of the Gromov boundary
$\partial\cT(\omega)$, as stated in \cref{it:tree-dims}.
\begin{restatement}{it:tree-dims}
    Let $Z_k$ be a Galton--Watson process with finitely supported offspring distribution with mean $m$ and maximal offspring number $N$.
    Let $\partial\cT$ denote the Gromov boundary of the associated Galton--Watson tree.
    Write
    \begin{equation*}
        \psi(R)=\frac{\log\log(1/R)}{\log(1/R)}.
    \end{equation*}
    The following results hold almost surely conditioned on non-extinction.

    For any dimension function $\phi$, if $\lim_{R\to 0}\frac{\psi(R)}{\phi(R)}=\alpha\in[0,\log N]$, then
    \begin{equation}\label{eq:first}
        \dimuAs\phi\partial\cT=\dimAs\phi \partial\cT = \alpha\left(1-\frac{\log m}{\log N}\right) +\log m.
    \end{equation}
    Otherwise, if $\lim_{R\to 0}\frac{\psi(R)}{\phi(R)}\geq\log N$, then
    \begin{equation}\label{eq:second}
        \dimA\partial\cT=\dimAs\phi\partial\cT =\log N.
    \end{equation}
\end{restatement}
\begin{proof}
    We may assume that the offspring number $X$ is not a constant almost surely, since otherwise $\log m = \log N$ and the theorem holds trivially.

    \smallskip
    Fix $\alpha \in (0,\log N)$.
    Recall that $\psi_\alpha(R)=\psi(R)/\alpha$, and that $\gamma$ is such that $m^\gamma=N$.
    We first show that
    \begin{equation}\label{eq:upperBd}
        \dimAs{\psi_\alpha} \partial \cT \leq \alpha\left( 1-\frac{\log m}{\log N} \right)+\log m\eqqcolon s_{\alpha}.
    \end{equation}
    Since $\alpha<\log N$, we can fix $\eps>0$ small enough that $t<\gamma$, where $t\coloneqq s_{\alpha}/\log m +2\eps$.

    First, let $E'_n$ be the event that a Galton--Watson tree has more than $m^{t n}$ descendants at level $n$, and write $E_k$ for the event that a Galton--Watson tree has more than $m^{tn}$ descendants at some level $n$ satisfying $\log (k-1) \leq \alpha n \leq \log k$.
    Note that the choice of $t$ guarantees that
    \begin{equation*}
        m^{(t-1-\eps)\frac{\gamma}{\gamma-1}\cdot\frac{\log(k-1)}{\alpha}}=(k-1)^{1+\eps\cdot\frac{\log m}{\alpha}\cdot\frac{\gamma}{\gamma-1}}.
    \end{equation*}
    Since $m^t<N$, by \cref{p:GW-large-deviations}, applying the above substitution yields
    \begin{align}
        P_k\coloneqq{}&\Prob(E_k)\nonumber\\
        \leq{}& \sum_{n=\lceil\frac{1}{\alpha}\log (k-1)\rceil}^{\lfloor \frac{1}{\alpha} \log k \rfloor}
    \Prob(E'_n)
    \label{eq:sum1}
    \\
        \lesssim{}&
        \frac{1}{\alpha}\left(\log k - \log (k-1)+1\right)\exp\left( -m^{(t-1-\eps)\frac{\gamma}{\gamma-1} \frac{\log(k-1)}{\alpha} } \right)
    \nonumber
    \\
        \lesssim{}& \exp\left( -(k-1)^{1+\eps\cdot\frac{\log m}{\alpha}\cdot\frac{\gamma}{\gamma-1}}\right)\nonumber,
    \end{align}
    with implicit constants independent of $k$.
    Thus since $\eps>0$,
    \begin{equation*}
      \sum_{k=1}^{\infty}P_k m^k \lesssim \sum_{k=1}^\infty \exp\left( k \log m - (k-1)^{1+\eps\cdot \frac{\log m}{\alpha}\cdot \frac{\gamma}{\gamma-1}} \right) < \infty.
    \end{equation*}
    Now by \cref{l:tree-bc}~\cref{im:BC0}, almost surely there exists a constant $C>0$ such that for all integers $k,n$ (with $k \geq 2$) satisfying $\log (k-1) \leq \alpha n \leq \log k$, every subtree at level $k$ has no more than $C m^{tn}$ descendants at level $k+n$.
    Now suppose $0<R<1$ and consider the ball $B(x,R)\subseteq\partial \cT$.
    By definition of the metric, the ball is the full subtree of a node $v\in L_k$, where $k=\lceil-\log R\rceil$.
    Thus by the definition of $\psi$,
    \begin{equation*}
      N_{R^{1+\psi_\alpha(R)}}(B(x,R))\lesssim m^{\frac{t \log k}{\alpha}}= k^{\frac{s_{\alpha}}{\alpha}+2\eps\frac{\log m}{\alpha}}\lesssim R^{-\psi_\alpha(R)(s_\alpha+2\eps\log m)}
    \end{equation*}
    with implicit constants independent of $R$.
    Since $\eps>0$ was arbitrary, the bound \cref{eq:upperBd} holds.

    \medskip
    We now give a proof of the lower bound.
    Let $\alpha \in(0,\log N)$, fix $\eps>0$, and let $t = s_{\alpha}/\log m-2\eps$.
    Taking $\eps$ to be sufficiently small guarantees that $t>0$.
    Let $E_k$ be the event that a Galton--Watson tree has at least $\tfrac{1-q}{2}m^{t(\log k)/\alpha}$
    descendants at level $\lfloor(\log k)/\alpha\rfloor$ that do not die out, recalling that $q$ is given in \cref{p:pgf}.
    By \cref{c:GW-deviations-surviving} and the choice of $t$, increasing the constant $C$ if necessary, we have
    \begin{equation*}
        P_k\coloneqq\Prob(E_k) \gtrsim \exp\left( -m^{(t-1+\eps)\frac{\gamma}{(\gamma-1)}\frac{\log k}{\alpha}} \right)
        = \exp\left( -k^{1-\eps\frac{\gamma\log m}{\alpha(\gamma-1)}} \right).
    \end{equation*}
    Thus,
    \begin{equation*}
        P_k m^k \gtrsim \exp\left( k\log m - k^{1-\eps\frac{\gamma\log m}{\alpha(\gamma-1)}} \right)
        \gtrsim \exp\left( \tfrac12 k \log m \right)
        = m^{k/2}.
    \end{equation*}
    Letting $K_k = m^{-k/2}$, we get
    \begin{equation*}
        \sum_{k=1}^{\infty}K_k P_k m^k\gtrsim \sum_{k=1}^{\infty}1 = \infty
        \quad
        \text{and}
        \quad
        \sum_{k=1}^\infty K_k = \sum_{k=1}^{\infty}m^{-k/2}< \infty.
    \end{equation*}
    By \cref{l:tree-bc}~\cref{im:BC1}, almost surely, conditioned on non-extinction, there exist infinitely many subtrees at levels $k_i\to\infty$ that have at least $\tfrac{1-q}{2}m^{t(\log k_i)/\alpha}$ descendants at level $k_i+\lfloor\tfrac{1}{\alpha}\log k_i\rfloor$ which do not die out.
    Let $R_i = e^{-k_i}$ and chose $x\in \partial \cT$ such that $B(x,R_i)$ is such a subtree at level $k_i$.
    Then since $R_i^{1+\psi_\alpha(R_i)}=e^{-k_i-\tfrac{1}{\alpha}\log k_i}$,
    \begin{equation*}
        N_{R_i^{1+\psi_\alpha(R_i)}}(B(x,R_i)) \geq\tfrac{1-q}{2}m^{t(\log k_i)/\alpha}
        =\tfrac{1-q}{2}
        k_i^{\frac{s_{\alpha}}{\alpha} - 2\eps\frac{\log m}{\alpha}}
        =\tfrac{1-q}{2}R_i^{-\psi_\alpha(R_i)( s_{\alpha} - 2\eps \log m)}.
    \end{equation*}
    Hence, as $\eps>0$ was arbitrary, almost surely $\dimAs{\psi_\alpha} \partial\cT \geq s_{\alpha}$.
    We note that the bound $\dimAs{\psi_\alpha}\partial\cT\geq\log N$ for $\alpha\geq\log N$ holds by a similar, but easier, argument.

    It remains to show that, almost surely, \cref{eq:first} and \cref{eq:second} hold simultaneously for any dimension function satisfying the hypotheses.
    First, fix a dense countable subset $\mathcal{C}=\Q\cap(0,\infty)$ so that, almost surely for all $\alpha\in \mathcal{C}$ simultaneously,
    \begin{equation*}
        \dimAs{\psi/\alpha}\partial\cT=\min\{\alpha,\log N\}\cdot\left(1-\frac{\log m}{\log N}\right)+\log m
    \end{equation*}
    and $\dimB\partial\cT=\log m$ and $\dimA\partial\cT=\log N$.

    Fix a typical element $\partial\cT$ as above.
    Note that, since $\dimAs{\psi_\alpha}\partial\cT$ is a continuous and increasing function of $\alpha$, by \cref{it:fund-dim}~\cref{im:limit-eq} and \cref{it:upper-recover}, the formula \cref{eq:first} holds for all dimension functions $\phi$ with $\lim_{R\to 0}\psi(R)/\phi(R)\in(0,\infty)$.
    Otherwise, if $\phi$ is a dimension function with $\lim_{R\to 0}\psi(R)/\phi(R)=0$, then for all $\alpha>0$,
    \begin{equation*}
        \dimB\partial\cT\leq\dimAs{\phi}\partial\cT\leq\dimuAs{\phi}\partial\cT\leq\alpha\left(1-\frac{\log m}{\log N}\right)+\log m,
    \end{equation*}
    and the upper bound converges to $\dimB\partial\cT$ as $\alpha$ converges to zero.
    Finally, suppose $\phi$ is a dimension function with $\lim_{R\to 0}\psi(R)/\phi(R)=\infty$.
    Then by \cref{it:upper-recover},
    \begin{equation*}
        \dimA\partial\cT\geq\dimAs{\phi}\partial\cT\geq\dimAs{\psi/(\log N)}\partial\cT=\dimA\partial\cT.
    \end{equation*}
    Thus the desired formulas hold for $\partial\cT$, as claimed.
\end{proof}
\begin{remark}
    In the proof above we assumed that the limit $\lim_{R\to 0}\tfrac{\psi(R)}{\phi(R)}=\alpha$ exists.
    For the lower bound, the existence of the limit was used in establishing the relationship between all (random) scales $r_i$ and $R_i$, which is a greater degree of independence than strictly necessary.
    For the upper bound, establishing summability of the $P_k m^k$ requires the range of the summation in \cref{eq:sum1} to be constrained appropriately, which in turn requires that the dimension function $\phi$ is close to the dimension function $\psi$ over large ranges of scales.
    Relaxing the assumptions in \cref{it:tree-dims} somewhat is possible, but obtaining a precise formula for all dimension functions will require substantially more work beyond what is done in this proof.
\end{remark}
\subsection{Mandelbrot percolation of the unit cube}
Galton--Watson processes are frequently used to model stochastically self-similar sets that arise from percolation processes.
One particularly notable example is that of Mandelbrot percolation of the unit cube in $\R^d$.
Fix an integer $n\geq 2$ and retention probability $p>n^{-d}$.
Let $M_0 = [0,1]^d$ be the unit cube and write $\mathbf{Q}'_1$ for the collection of $n^d$ subcubes of $M_0$ of side-length $1/n$ that evenly partition $M_0$.
For each subcube $Q\in\mathbf{Q}'_1$ we independently decide to keep it with probability $p$.
Call this collection of subcubes $\mathbf{Q}_1$ and set $M_1 = \bigcup\mathbf{Q}_1 \subseteq M_0$.
Having constructed $\mathbf{Q}'_1$, $\mathbf{Q}_1$ and $M_1$, we iteratively construct $\mathbf{Q}'_{k+1}$ to be the set of all $n^d  \cdot \#\mathbf{Q}_k$ subcubes of the cubes in $\mathbf{Q}_k$, for $k \geq 1$.
The collection $\mathbf{Q}_{k+1}$ is the set of independently retained cubes in $\mathbf{Q}'_{k+1}$ with probability $p$.
Finally, we set $M_{k+1}=\bigcup \mathbf{Q}_{k+1}'$.
The process is known as Mandelbrot percolation and the limit set $M = \bigcap_{k\in\N} M_k \subset \R^d$ is a stochastically self-similar set, see \cref{fig:Mandelbrot}.
Note that the number of subcubes of a cube $Q\in\mathbf{Q}_k$ has binomial distribution $B(n^d,p)$ and is independent of other subcubes.
Thus, $\#\mathbf{Q}_k$ is a Galton--Watson process with offspring distribution $B(n^d,p)$ and we may index the cubes using a Galton--Watson tree $\partial T$.

Mandelbrot percolation is a special case of general fractal percolation.
We may consider the unit cube as the invariant set under the IFS of similarities $S_\tau(x) = n^{-1}( x + t_{\tau})$, where $t_\tau\in\left\{ 0,1,\ldots,n-1 \right\}^d$.
In general, one may consider such percolation on any attractor given by an IFS, see \cite{zbl:1185.28013} which first explicitly considered such construction.
Given such strong correspondence with Galton--Watson processes, we may immediately apply \cref{it:tree-dims} to all fractal percolation whose cylinder sets are of controlled size such as homogeneous self-similar iterated function systems satisfying the open set condition.
Using slightly more general considerations, they also apply to non-homogeneous systems, see e.g.~\cite{zbl:1437.28015} for a detailed discussion.
To keep proofs succinct we will only prove the case for Mandelbrot percolation in \cref{ic:percolation}.
\begin{proofref}{ic:percolation}
    Consider a Galton--Watson process $Z_k$ with binomial offspring distribution $X\equiv_D B(n^d, p)$ and its associated Galton--Watson tree $\cT$ with Gromov boundary $\partial\cT$.
    It is straightforward to verify that $m=\E(X) = p n^d$ and $N=n^d$.
    Further, changing the metric $d(x,y)=e^{|x\wedge y|}$ on $\partial\cT$ to $d'(x,y)=d(x,y)^{\log n}=n^{|x\wedge y|}$ has the effect of changing the conclusion of \cref{it:tree-dims} to the following.

    If $\lim_{R\to 0}\frac{\psi(R)}{\phi(R)}=\alpha\in[0,\log N] = [0,\log n^d]$,
    \begin{equation*}
        \dimAs{\phi} \partial T = \frac{1}{\log n}\left(\alpha \left( 1-\frac{\log m}{\log N}
        \right)+\log m\right)
        =\alpha \frac{\log (1/p)}{d\log^2 n}+ \frac{\log p n^d}{\log n}.
    \end{equation*}
    Otherwise if $\liminf_{R\to 0}\frac{\psi(R)}{\phi(R)}\geq \log n^d$,
    \begin{equation*}
        \dimAs\phi \partial T = \log_n N = d.
    \end{equation*}

    To see that the same conclusion may also be reached for Mandelbrot percolation, note that the Galton--Watson process $Z_k$ counts the number of surviving subcubes at iteration level $k$.
    In particular, there exists a natural bijection between vertices $v\in\cT$ and surviving subcubes $Q_v$ of side-length $n^{-|v|}$ such that if $w=vj\in\cT$ for some $j\in\{1,\ldots,N\}^*$ then $Q_{w}\subseteq Q_v$ is a subcube of side-length $n^{-|w|}$.
    Clearly $M = \cap_{k\in\N}\cup_{v\in L_k}Q_v$.

    Write $Q(x,k)$ for a subcube of side-length $n^{-k}$ containing $x$.
    Then for all $x\in M$ and $0<r<1$,
    \begin{equation*}
        B(x,r) \supseteq Q(x,\lceil \log_n (\sqrt{d}/r) \rceil)
    \end{equation*}
    and $B(x,r) \subseteq \bigcup_{j=1}^{3^d}Q_{v_j}$ where $Q_{v_1}=Q(x,\lfloor \log_n (1/r)\rfloor)$ and $Q_{v_j}$ are the (at most) $3^d-1$ neighbouring subcubes of level $\lfloor \log_n (1/r)\rfloor$.
    Note further that the subtree $\cT(v)$ that corresponds to $Q_v$ is of diameter $n^{-|v|}$ making their corresponding diameters comparable:
    \begin{equation*}
        2r=\diam(B(x,r))\approx n^{\log_n r} \approx n^{-\lceil \log_n (\sqrt{d}/r) \rceil} =\diam(\cT(v)).
    \end{equation*}
    Therefore
    \begin{equation*}
      N_r(\cT(w)) \lesssim N_r(B(x,R)) \lesssim \sum_{j=1}^{3^d} N_{r}(\cT(v_j)),
    \end{equation*}
    where $Q_{v_1} = Q(x,\lfloor \log_n(1/r)\rfloor)$ and $Q_w = Q(x,\lceil \log_n (\sqrt{d}/r) \rceil)$.
    The corollary now follows directly from \cref{it:tree-dims}.
\end{proofref}

\section{Self-similar sets and decreasing sequences}\label{s:other-ex}
\subsection{Self-similar sets and general upper bounds}\label{ss:similar-upper}
Let $\mathcal{I}$ be a finite index set and let $\{S_i\}_{i\in\mathcal{I}}$ be an iterated function system (IFS) of similarities, i.e.
\begin{equation*}
    S_i(x)=r_i x+d_i\text{ for }0<|r_i|<1\text{ and }d_i\in\R.
\end{equation*}
Let $\mathcal{I}^*=\bigcup_{n=0}^\infty\mathcal{I}^n$ and for $r>0$, let
\begin{equation*}
    \Lambda_r=\{ \sigma\in\mathcal{I}^*:r_\sigma\leq r<r_{\sigma^-} \},
\end{equation*}
where if $\sigma = \sigma_1 \dotsb \sigma_n$ then $r_{\sigma} = r_{\sigma_1} \dotsb r_{\sigma_n}$ and $\sigma^- = \sigma_1 \dotsb \sigma_{n-1}$.
Set
\begin{equation}\label{e:Mr-def}
    \mathcal{M}_r(x)=\{ S_\sigma:\sigma\in\Lambda_r,S_\sigma(K)\cap B(x,r)\neq\varnothing \}\qquad\text{and}\qquad M_r=\sup_{x\in K}\#\mathcal{M}_r(x).
\end{equation}
This is a generalization of the notation $\widetilde{M}_n$ from the introduction required to handle non-homogeneous iterated function systems.
We then say that the IFS satisfies the \defn{weak separation condition} (WSC) if $\sup_{r>0}M_r<\infty$.\footnote{
    The weak separation condition was first introduced in \cite{zbl:0929.28007} with a somewhat different definition; this version was proven to be equivalent in \cite{zbl:0874.54025}.
}
In this case, the attractor $K$ is Ahlfors--David regular and $\dimH K=\dimA K$.
Otherwise, $\dimA K=1$ \cite[Theorem~1.3]{zbl:1317.28014}.

It is straightforward to use the counts $M_r$ to give a general upper bound relevant for any self-similar IFS.
In the special case that the WSC condition holds, $M_r$ is uniformly bounded above and this result just shows that $\dimAs\phi K=\dimH K$ for all dimension functions $\phi$.
\begin{proposition}\label{p:self-similar-upper}
    Let $\phi$ be any function such that
    \begin{equation*}
        \lim_{r\to 0}\frac{\log M_r}{\phi(r)\log(1/r)}=0.
    \end{equation*}
    Then $\dimAs\phi K=\dimH K$.
\end{proposition}
\begin{proof}
    Let $0<r<1$ and let $x\in K$.
    First, recalling \cite[Theorem~4]{zbl:0683.58034}, there is a constant $C>0$ (not depending on $r$) so that with $N = Cr^{-\phi(r)\dimB K}$, there is a cover $\{B(y_i,r^{\phi(r)})\}_{i=1}^N$ for $K$.
    Therefore by self-similarity,
    \begin{equation*}
        \{S_\sigma(B(y_i,r^{\phi(r)}))\}_i=\{B(S_\sigma(y_i),r_\sigma r^{\phi(r)})\}_i
    \end{equation*}
    is a cover for $S_\sigma(K)$ with balls of radius $r_\sigma r^{\phi(r)}$.
    Thus since $K \cap B(x,r)\subset\bigcup_{S_\sigma\in\mathcal{M}_r(x)}S_\sigma(K)$, applying the above observation to each image $S_\sigma(K)$,
    \begin{equation*}
        N_{r^{1+\phi(r)}}(K \cap B(x,r))\leq C M_r\left(\frac{r}{r^{1+\phi(r)}}\right)^{\dimB K}.
    \end{equation*}
    Moreover, for every $\eps>0$ and all $r$ sufficiently small (depending on $\eps$), by assumption on $\phi$,
    \begin{equation*}
        M_r\leq\left(\frac{r}{r^{1+\phi(r)}}\right)^\eps.
    \end{equation*}
    It follows that $\dimAs\phi K=\dimB K=\dimH K$.
\end{proof}

\subsection{Lower bounds and Assouad dichotomy}
A general strategy for obtaining lower bounds for the Assouad dimension of a set is to construct subsets which are, in some sense, close to being arithmetic progressions.
\begin{definition}
    Let $K\subset\R$ be a compact set.
    For $r,\eps>0$ we call a set $(r^{-1} K -t)\cap [0,1]$ which is an $\eps$-dense subset of $[0,1]$ an \emph{$(r,\eps)$-microset} of $K$.
\end{definition}
One can think of the parameter $r$ as the \emph{scale}, and the parameter $\eps$ as the \emph{resolution}.
Of course, if $K$ has an $(r,\eps)$-microset and $c\in(0,1)$, then $K$ also has a $(cr, c^{-1}\eps)$-microset.
Moreover, $K$ has Assouad dimension 1 if and only if $K$ has an $(r,\eps)$-microset for arbitrarily small $\eps$ \cite{doi:10.1093/imrn/rnw336}.

It is proven in \cite[Theorem~1.3]{zbl:1317.28014} that any self-similar set in $\R$ that is not a singleton either satisfies the WSC, or has Assouad dimension 1.
In fact, we can interpret the proof as a certain ``sequential amplification'' of microsets.
Let $\{S_i\}_{i\in\mathcal{I}}$ be an IFS of similarities with attractor $K$ such that $S_i(x)=\rho x + d_i$ for $i\in\mathcal{I}$ where $\rho\in(0,1)$ is fixed.
Define a distance on $\mathcal{I}^*$ by
\begin{equation*}
    d(\sigma,\tau)=\begin{cases}
      \rho^{-n}|S_\sigma(0)-S_\tau(0)| &\text{if }|\sigma|=|\tau|=n,\\*
      \infty &\text{otherwise.}
    \end{cases}
\end{equation*}
Here, $|\sigma |$ denotes the length of the word $\sigma$.
We begin with a simple lemma demonstrating this amplification process for equicontractive self-similar sets.
\begin{lemma}\label{l:tan-amplify}
    Let $\{S_i\}_{i\in\mathcal{I}}$ be an equicontractive IFS of similarities with attractor $K$.
    Suppose $K$ has an $(r,\eps)$-microset.
    Then if $n\in\N$ and $c>0$ are such that there are words $\sigma,\tau\in\mathcal{I}^n$ such that $c r<d(\sigma,\tau)\leq r$, then $K$ contains a $\bigl((1+c)\rho^n r, (1+c)^{-1}\eps\bigr)$-microset.
\end{lemma}
\begin{proof}
    Let $t\in\R$ and $T(x)=r^{-1} x-t$ be a similarity so that with $P\coloneqq T^{-1}([0,1])\cap K$, $T(P)$ is an $(r,\eps)$-microset.
    Next, by assumption, get $n\in\N$ and words $\sigma,\tau\in\mathcal{I}^n$ such that $c<\delta\leq 1$ where
    \begin{equation*}
        \delta\coloneqq r^{-1}\rho^{-n}(S_\tau(0)-S_\sigma(0)),
    \end{equation*}
    Write $Q=S_\sigma(P)\cup S_\tau(P)$.
    Of course, $Q\subset K$ by self-similarity and moreover a direct computation gives that
    \begin{equation*}
        T\circ S_\sigma^{-1}(Q) = T(P)\cup (T(P)+\delta).
    \end{equation*}
    Since $T(P)$ is an $\eps$-dense subset of $[0,1]$ and $\delta\leq 1$, it follows that $\bigl(T(P)\cup(T(P)+\delta)\bigr)\cap[0,1+\delta]$ is an $\eps$-dense subset of $[0,1+\delta]$.
    Therefore writing $h(x)=(1+c)^{-1} x$ and recalling that $\delta>c$,
    \begin{equation*}
        h\circ T\circ S_\sigma^{-1}(K)\cap[0,1]
    \end{equation*}
    is an $\eps/(1+c)$-dense subset of $[0,1]$.
    Moreover, $h\circ T\circ S_\sigma$ has similarity ratio $(1+c)^{-1}r^{-1}\rho^{-n}$, which gives the claim.
\end{proof}
Of course, having precise information about the existence of microsets gives lower bounds for the $\phi$-Assouad dimensions.
\begin{lemma}\label{l:micro-bounds}
    Suppose $K$ has an $(r,\eps)$-microset.
    Then there exists $x\in K$ such that
    \begin{equation}\label{e:large-count}
        N_{\eps r}\bigl(K \cap B(x,r)\bigr)\geq\frac{1}{\eps}.
    \end{equation}
    In particular, if $K$ has a sequence of $(r_k,\eps_k)$-microsets with $\eps_k$ converging to zero and $\phi$ is any dimension function such that
    \begin{equation}\label{e:phik-control}
        \limsup_{k\to\infty}\frac{\log\eps_k}{\phi(r_k)\log r_k}\geq 1,
    \end{equation}
    then $\dimuAs\phi K=1$.
\end{lemma}
\begin{proof}
    Note that \cref{e:large-count} follows directly from the definition of an $(r,\eps)$-microset.
    Now re-arranging \cref{e:phik-control} gives that
    \begin{equation*}
        r_k^{1+\phi(r_k)}\geq\eps_k r_k
    \end{equation*}
    or equivalently by \cref{e:large-count} there is an $x\in K$ so that
    \begin{equation*}
        N_{\eps_k r_k}\bigl(B(x,r_k) \cap K) \geq\frac{1}{\eps_k}.
    \end{equation*}
    But this holds for infinitely many $k$ with $\eps_k$ converging to $0$, so $\dimuAs\phi K=1$.
\end{proof}

Suppose the WSC fails, which implies that $0$ is an accumulation point of the set $\{d(\sigma,\tau):\sigma,\tau\in\mathcal{I}^*\}\setminus\{0\}$ \cite{zbl:0874.54025}.
Moreover suppose, inductively, that we have constructed an $(r,\eps)$-microset.
Intending to apply \cref{l:tan-amplify}, get $n_0$ and $\sigma,\tau\in\mathcal{I}^{n_0}$ such that $d(\sigma,\tau)\leq r$.
Of course, $d(\sigma,\tau)$ may be much smaller than $r$; however, if $i_0\in\mathcal{I}$ is any fixed letter, then $d(\sigma i_0,\tau i_0)=\rho^{-1}d(\sigma,\tau)$.
Repeatedly appending the letter $i_0$ guarantees the existence of some $n\geq n_0$ and words $\sigma,\tau\in\mathcal{I}^n$ satisfying the hypotheses of \cref{l:tan-amplify} with $c=\rho$, so $K$ in fact contains a $((1+\rho) r \rho^n,(1+\rho)^{-1}\eps)$-microset.
But $(1+\rho)^{-1}<1$ is a fixed constant, so repeating this construction we see that $K$ has an $(r,\eps)$-microset for arbitrarily small $\eps$, so $\dimA K=1$ under the assumption that the WSC fails\footnote{This gives a somewhat more transparent proof of \cite[Theorem~3.1]{zbl:1317.28014} and \cite[Theorem~4.1]{zbl:1441.28003} in the self-similar equicontractive case.}.

In order to understand the rate at which the $(r,\eps)$-microsets converge to $[0,1]$ in the Hausdorff metric, we would like to make the above proof quantitative.
However, the main challenge is to control the level $n\geq n_0$ at which we have a gap of precisely the correct size.
Thus it seems to be challenging to show in general that the upper bound in \cref{p:self-similar-upper} is sharp.
In the next section, we construct explicit examples of self-similar sets for which the gap $n\geq n_0$ can be controlled very precisely.
Using this we can prove relatively good bounds on the resolution at which large microsets appear.

\subsection{Explicit lower bounds for a class of self-similar sets}
Given some $m\in\N$ with $m\geq 3$ and $t>0$, consider the homogeneous self-similar IFS defined by
\begin{equation}\label{e:st-def}
    S_0(x)=\frac{x}{m}\qquad S_t(x)=\frac{x}{m}+t\qquad S_1(x)=\frac{x+1}{m}.
\end{equation}
Denote the attractor of this IFS by $K$.
Note that when $t$ is irrational, this IFS satisfies the exponential separation condition (see \cite[Theorem~1.6]{zbl:1337.28015} for an explicit statement and proof) and $\dimH K=\tfrac{\log 3}{\log m}$.
As described in the proof of \cite[Theorem~7.3.1]{zbl:1467.28001}, one can use work of Shmerkin \cite{zbl:1426.11079} together with the fact that the IFS satisfies the exponential separation condition to show that $\dimqA K = \dimH K$.
Observe that $C\subset K$, where $C$ is the Cantor set consisting of the points
\begin{equation*}
    C=\left\{\sum_{n=1}^\infty\frac{j_n}{m^n}:j_n\in\{0,1\}\right\}.
\end{equation*}
Note that $C$ does not depend on the choice of $t$.
Recalling that $S_0(0)=0$ so that $0\in K$, we denote the set of endpoints at level $k$ by
\begin{equation*}
    E_k\coloneqq\{ S_\sigma(0):\sigma\in\{0,t,1\}^k \}\subset K.
\end{equation*}
Of course, $\{0\}=E_0\subset E_1\subset E_2\subset\cdots$.

We now introduce the key definition relevant for our construction.
\begin{definition}\label{n:t-def}
    Suppose we have chosen a sequence of positive integers $(n_k)_{k=1}^\infty$.
    Write $N_k=n_1+\cdots+n_k$.
    For each $k\in\N$, let $\omega_k\in\{0,1\}^{n_k}$ denote the word $\omega_k=(0,\ldots,0,1)$ where $0$ is repeated $n_k-1$ times.
    We then let $t = t(n_k)_{k=1}^\infty \in (0,1/(m-1)]$ denote the point with base-$m$ expansion $\omega_1\omega_2\ldots$ and let $K = K(n_k)_{k=1}^\infty$ denote the attractor of the IFS $\{S_0,S_t,S_1\}$.
\end{definition}
We also recall that $C$ is the associated Cantor set, which does not depend on the choice $t$.
The main work in this section is to choose the sequence $(n_k)_{k=1}^\infty$ so that $K$ has large microsets at precise scales.
\begin{lemma}\label{l:nk-choice}
    Let $m\in\N$ with $m\geq 3$ be fixed.
    Then there exists a sequence $(n_k)_{k=1}^\infty$ such that the corresponding attractor $K(n_k)_{k=1}^\infty$ has an $\bigl(r_k,(1+1/m^2)^{-k}\bigr)$-microset for each $k\in\N$, where for each $\eps>0$,
    \begin{equation*}
        -\log r_k\leq\frac{2}{\eps}\cdot(2+\eps)^{k-1}\log m.
    \end{equation*}
\end{lemma}
\begin{proof}
    We first define a sequence $(n_k)_{k=1}^\infty$ of positive integers and a sequence $(r_k)_{k=1}^\infty$ of real numbers inductively.
    Begin with $n_1=1$ and $r_1=1$.
    Now suppose we have chosen $n_1,\ldots,n_k$ and $r_1,\ldots,r_k$ for some $k\in\N$.
    Let $n_{k+1}$ be chosen so that
    \begin{equation}\label{e:nk-choice}
        \frac{r_k}{m}< \frac{m}{m^{n_{k+1}}}\leq r_k
    \end{equation}
    and define $r_{k+1}=(1+1/m^2)r_k m^{-N_k}$.

    Now let $N_k$, $\omega_k$, and $t$ be as in \cref{n:t-def} and let $K$ denote the corresponding attractor.
    Note that $(r_k)$ is strictly decreasing and $(n_k)$ is strictly increasing, so $t$ is irrational.

    For each $k\in\N\cup\{0\}$, let $z_k\in E_{n_k}\cap C$ be the point coded by the word $\omega_1\omega_2\dotsb \omega_k\in\mathcal{I}^{N_k}$.
    Then $m^{N_k}(t-z_k)$ has base-$m$ expansion $\omega_{k+1}\omega_{k+2}\dotsb$.
    In particular, the words $\sigma_k=\omega_1\dotsb\omega_k$ and $\tau_k=t0\dotsb 0$ in $\mathcal{I}^{N_k}$ satisfy
    \begin{equation}\label{e:nk-param}
        1<m^{n_{k+1}}d(\sigma_k,\tau_k)=1+\frac{1}{m^{n_{k+2}}}+\frac{1}{m^{n_{k+2}+n_{k+3}}}+\cdots\leq m.
    \end{equation}

    Next, we verify that $K$ contains an $(r_k,(1+1/m^2)^{-k})$-microset for all $k\in\N$.
    Of course, the $k=1$ case follows immediately since $\{0,1/(m-1)\}\subset K$.
    Now let $k\in\N$ be arbitrary and suppose $K$ contains an $(r_k,(1+1/m^2)^{-k})$-microset.
    Combining \cref{e:nk-param} and \cref{e:nk-choice}, there are words $\sigma_k$ and $\tau_k$ in $\mathcal{I}^{N_k}$ so that
    \begin{equation*}
        \frac{r_k}{m^2}< d(\sigma_k, \tau_k)\leq r_k.
    \end{equation*}
    Therefore applying \cref{l:tan-amplify} and the definition of $r_{k+1}$, $K$ contains a $\bigl(r_{k+1},(1+1/m^2)^{-(k+1)}\bigr)$-microset, as required.

    Finally, we must lower bound $r_k$.
    For notational simplicity, write $c(\eps)=\frac{2}{\eps(2+\eps)}$.
    We will prove by induction that for all $\eps>0$ and $k\in\N$ that
    \begin{equation*}
        N_k \leq c(\eps)(2+\eps)^k\qquad\text{and}\qquad -\log r_k\leq c(\eps)(2+\eps)^k\log m.
    \end{equation*}
    First, note that $c(\eps)$ is chosen precisely so that $2c(\eps)(2+\eps)^k+2\leq c(\eps)(2+\eps)^{k+1}$ for all $k \in \N$.
    Now the case $k=1$ holds since $n_1=1$ and $-\log r_1=0$.
    Suppose the hypothesis holds for $k\in\N$.
    Then directly from the definitions,
    \begin{align*}
        N_k +n_{k+1}&\leq N_k +\frac{\log\left(\frac{m^2}{r_k}\right)}{\log m}\\
                              &\leq c(\eps)(2+\eps)^k+c(\eps)(2+\eps)^k+2\\
                              &\leq c(\eps)(2+\eps)^{k+1}.
    \end{align*}
    Moreover,
    \begin{align*}
        -\log r_{k+1} &= N_k\log m -\log r_k-\log\left(1+\frac{1}{m^2}\right)\\
                      &\leq c(\eps)(2+\eps)^k\log m +c(\eps)(2+\eps)^k\log m \\
                      &\leq c(\eps)(2+\eps)^{k+1}\log m,
    \end{align*}
    as claimed.
\end{proof}
Recalling the definition of $M_r$ for $r>0$ from \cref{e:Mr-def}, define ${\widetilde M}_k=M_{m^{-k}}$ for non-negative integers $k$.
This is the same as the definition from the introduction.
\begin{theorem}\label{t:ss-lower}
    For all $m\in\N$ with $m \geq 3$ there exists $t\in(0,1)$ such that the attractor $K$ of the IFS $\{S_0,S_t,S_1\}$ defined in \cref{e:st-def} has the following properties.
    If $\psi$ is a dimension function satisfying
    \begin{equation*}
        \lim_{n\to\infty}\frac{\log\widetilde{M}_n}{n\cdot\psi(m^{-n})}=0,
    \end{equation*}
    then $\dimuAs\phi K=\dimH K$, and if $\phi$ is a dimension function satisfying
    \begin{equation}\label{e:limsup-n}
        \liminf_{n\to\infty}\frac{\log n}{n\cdot \phi\bigl(m^{-n}\bigr)}\geq \frac{(\log 3)\cdot(\log m)}{\log(1+1/m^2)},
    \end{equation}
    then $\dimAs\phi K=\dimA K=1$.
\end{theorem}
\begin{proof}
    The first claim is just a restatement of \cref{p:self-similar-upper} in slightly different notation.

    We will now verify the second claim.
    Applying \cref{l:nk-choice}, get a sequence $(n_k)_{k=1}^\infty$ and corresponding value $t$ as in \cref{n:t-def} so that the attractor $K$ has an $(r_k, (1+1/m^2)^{-k})$-microset for all $k\in\N$ where (taking $\eps=1$)
    \begin{equation*}
        -\log r_k\leq 2\cdot 3^{k-1}\log m.
    \end{equation*}
    For each $k\in\N$, let $\ell_k$ be such that $r_k / m <m^{-\ell_k} \leq r_k$.
    Then $K$ has an $(m^{-\ell_k}, m(1+1/m^2)^{-k})$-microset for all $k\in\N$ and moreover for all $k$ sufficiently large,
    \begin{equation}\label{e:lk-bound}
        \ell_k \leq 3^{k}.
    \end{equation}
    Next, let $\phi$ be a dimension function satisfying \cref{e:limsup-n}.
    Recalling \cref{im:phi-gap} in the definition of a dimension function (see \cref{d:dimf}),
    \begin{equation*}
        \lim_{k\to\infty}\frac{1}{\ell_k \phi(m^{-\ell_k})} = 0.
    \end{equation*}
    Using this fact followed by \cref{e:limsup-n} and then \cref{e:lk-bound},
    \begin{align*}
        \limsup_{k\to\infty}\frac{\log\bigl(m(1+1/m^2)^{-k}\bigr)}{\phi\bigl(m^{-\ell_k}\bigr) \cdot \log\bigl(m^{-\ell_k}\bigr)}
        &\geq \limsup_{k\to\infty}\frac{k\log (1+1/m^2)}{(\ell_k \log m)\cdot \phi(m^{-\ell_k})}\\
        &\geq \limsup_{k\to\infty}\frac{k\log 3}{\log \ell_k}\\
        &\geq 1.
    \end{align*}
    Therefore by \cref{l:micro-bounds}, $\dimAs\phi K=\dimA K=1$, as required.
\end{proof}

\subsection{Decreasing sequences with decreasing gaps}\label{ss:decr-seq}
The following result describes the $\phi$-Assouad dimensions of decreasing sequences with decreasing gaps.
For any such sequence $F$ we have $\{ \dimqA F, \dimA F \} \subseteq \{0,1\}$ (see~\cite{zbl:1479.28006} and~\cite{zbl:1405.28007}).
In fact, the quasi-Assouad dimension is $0$ if and only if the upper box dimension is $0$, and the Assouad dimension is $0$ if and only if the sequence is lacunary.
\begin{definition}
    We say that a function $f\colon(1,\infty)\to(0,1)$ has \defn{regular gaps} if $f$ is differentiable and for all sufficiently large $x$:
    \begin{enumerate}[nl]
        \item $f(x)$ is strictly decreasing and converges to $0$,
        \item $f'(x)$ is strictly increasing and converges to $0$, and
        \item $f'(x+1)/f'(x)$ is increasing and converges to $1$.
    \end{enumerate}
\end{definition}
Given a function $f$ with regular gaps, we define an associated sequence set
\begin{equation*}
    F = F_f\coloneqq\{ f(n):n\in\N\}\cup\{0\}.
\end{equation*}
In \cref{t:decreasingseq} below, we establish a dimension result for the sets $F_f$.
\begin{remark}
    Examples of functions to which we can apply \cref{t:decreasingseq} include $f(x)=e^{-x^{\alpha}}$ for $0<\alpha<1$ and $f(x)=e^{-(\log x)^{\beta}}$ for $\beta>1$.
    The specific example $n^{-\log n}$ has been considered in \cite[Example~2.18]{zbl:1485.28006}.
    However, even for examples like these, if there is no simple closed form expression for $(f')^{-1}$ it may be difficult to get a more explicit formula for the $\phi$-Assouad dimension than \cref{e:dedgformula}.
\end{remark}
\begin{theorem}\label{t:decreasingseq}
    Let $f$ be a function with regular gaps and assume that the associated sequence set $F$ has upper box dimension $0$.
    Let $\phi$ be a dimension function with $\phi(R) \to 0$ as $R \to 0$.
    Write
    \begin{equation*}
        \beta=\limsup_{x\to\infty}\bigl((f(x))^{1+\phi(f(x))}+f'(x)\bigr).
    \end{equation*}
    If $\beta\geq 0$, then $\dimAs\phi F=1$.
    Otherwise,
    \begin{equation}\label{e:dedgformula}
        \dimAs{\phi} F =\limsup_{R \to 0^+} \frac{\log \left( R^{-1-\phi(R)} f\bigl((f')^{-1}(- R^{1+\phi(R)})\bigr)  + (f')^{-1}\bigl(- R^{1+\phi(R)}\bigr)  - R\right)}{-\phi(R) \log R}.
    \end{equation}
\end{theorem}

\begin{proof}
    The strategy to obtain~\cref{e:dedgformula} is to fix scales $0 < r < R$ and cover a whole interval around~$0$ with intervals of size $r$ until the point where the gap size exceeds $r$, and then cover each of the remaining points in $F \cap (0,R)$ individually.
    By the mean value theorem, $f(n) - f(n+1) = - f'(x_n)$ for some $x_n \in (n,n+1)$, so the sequence of gaps is strictly decreasing.
    Moreover, since $f'$ is monotonic and $f'(x+1)/f'(x) \to 1$ as $x \to \infty$,
    \begin{equation}\label{e:fn-lim}
        \lim_{n \to \infty} \frac{f(n) - f(n+1)}{- f'(n)} = 1.
    \end{equation}
    Since the sequence of gaps is decreasing,
    \begin{equation*}
        N_r((0,R) \cap F) \gtrsim N_r((x,x+R) \cap F)
    \end{equation*}
    uniformly for all $x \in \R$ and $0<r<R$, so it suffices to consider intervals whose left endpoint is~$0$.

    If $\beta\geq 0$, then for infinitely many $n$, writing $R_n = f(x_n)$, by \cref{e:fn-lim} and the definition of $\beta$
    \begin{equation*}
        f(\lfloor x_n \rfloor) - f(\lfloor x_n\rfloor+1) \lesssim-f'(\lfloor x_n\rfloor)\lesssim R_n^{1+\phi(R_n)}.
    \end{equation*}
    Since the gaps are decreasing, it follows that
    \begin{equation*}
        N_{R_n^{1+\phi(R_n)}}((0,R_n) \cap F) \gtrsim R_n^{-\phi(R_n)},
    \end{equation*}
    so $\dimAs{\phi} F = 1$.

    Otherwise, $\beta <0$.
    Recalling the formula \cref{e:phi-dimA-lim} for the $\phi$-Assouad dimension, it suffices to show that the expression inside the limit supremum in \cref{e:dedgformula} is asymptotically the same as expression inside the limit supremum in \cref{e:phi-dimA-lim}.
    Fix some small $R = f(x)$ and let $r \coloneqq R^{1+\phi(R)}$.
    Let $n \coloneqq \lfloor x \rfloor$ and $m \coloneqq \lceil (f')^{-1}(-r) \rceil$, so $R \approx f(n)$ and $r \approx f(m) - f(m+1)$.
    Note that the assumption $- f'(x) > (f(x))^{1+\phi(f(x))} = r$ from $\beta < 0$ means that $m > n$.
    Now,
    \begin{equation*}
        N_r ((0,R) \cap F) \approx \frac{f(m)}{r} + m - n.
    \end{equation*}
    Taking logarithms, \cref{e:dedgformula} follows from \cref{e:phi-dimA-lim}.
\end{proof}

\begin{acknowledgements}
AB was supported by an EPSRC New Investigators Award (EP/W003880/1) at Loughborough University.
Part of this work was completed while AB was a PhD student at the University of St Andrews, supported by a Leverhulme Trust Research Project Grant (RPG-2019-034).

AR was supported by EPSRC Grant EP/V520123/1 and the Natural Sciences and Engineering Research Council of Canada.
Part of this work was completed while AR visited the University of Oulu while supported by a scholarship from the London Math Society and the Cecil King Foundation.

ST was supported by the European Research Council Marie Skłodowska--Curie Personal Fellowship \#101064701.

The authors thank Kenneth Falconer and Jonathan Fraser for some helpful comments on a draft version of this document.
We also thank two anonymous referees for a thorough reading and many helpful comments.
\end{acknowledgements}
\end{document}